\documentclass[12pt,oneside]{article}%
\usepackage[centertags]{amsmath}
\usepackage{amssymb}
\usepackage{amsthm}
\usepackage{amsfonts}
\usepackage{graphicx}%
\usepackage{hyperref}
\providecommand{\U}[1]{\protect\rule{.1in}{.1in}}
\newtheorem{theorem}{Theorem}[section]
\newtheorem{corollary}{Corollary}[section]
\newtheorem{remark}{Remark}[section]
\newtheorem{example}{Example}[section]

\newtheorem{proposition}{Proposition}[section]
\newtheorem{definition}{Definition}[section]

\newcommand{\punt}{\boldsymbol{ . }}
\begin{document}

\title{The classical umbral calculus: Sheffer sequences}
\author{E. Di Nardo \thanks{Dipartimento di Matematica e Informatica, Universit\`a
degli Studi della Basilicata, Viale dell'Ateneo Lucano 10, 85100 Potenza,
Italia, elvira.dinardo@unibas.it}, H. Niederhausen \thanks{Department of
Mathematical Sciences, Florida Atlantic University, 777 Glades Road, Boca
Raton , Florida 33431-0991, USA, niederha@fau.edu }, D. Senato
\thanks{Dipartimento di Matematica e Informatica, Universit\`a degli Studi
della Basilicata, Viale dell'Ateneo Lucano 10, 85100 Potenza, Italia,
domenico.senato@unibas.it}}
\date{\today}
\maketitle

\begin{abstract}
Following the approach of Rota and Taylor \cite{SIAM}, we present an
innovative theory of Sheffer sequences in which the main properties are
encoded by using umbrae. This syntax allows us noteworthy computational
simplifications and conceptual clarifications in many results involving
Sheffer sequences. To give an indication of the effectiveness of the theory,
we describe applications to the well-known connection constants problem, to
Lagrange inversion formula and to solving some recurrence relations.

\end{abstract}

\textsf{\textbf{keywords}: umbral calculus, Sheffer sequences, connection
constants problem, linear recurrences, Lagrange inversion formula. }%

\section{Introduction}

As well known, many polynomial sequences like Laguerre polynomials, first and
second kind Meixner polynomials, Poisson-Charlier polynomials and Stirling
polynomials are Sheffer sequences. Sheffer sequences can be considered the
core of umbral calculus: a set of tricks extensively used by mathematicians at
the beginning of the twentieth century.

Umbral calculus was formalized in the language of the linear operators by
Gian-Carlo Rota in a series of papers (see \cite{Bellnumbers}, \cite{VIII},
and \cite{Roman}) that have produced a plenty of applications (see
\cite{dibucchianico}). In 1994 Rota and Taylor \cite{SIAM} came back to
foundation of umbral calculus with the aim to restore, in an light formal
setting, the computational power of the original tools, heuristical applied by
founders Blissard, Cayley and Sylvester. In this new setting, to which we
refer as the classical umbral calculus, there are two basic devices. The first
one is to represent a unital sequence of numbers by a symbol $\alpha,$ called
an umbra, that is, to represent the sequence $1, a_{1}, a_{2}, \ldots$ by
means of the sequence $1, \alpha, \alpha^{2}, \ldots$ of powers of $\alpha$
via an operator $E,$ resembling the expectation operator of random variables.
The second device is that distinct umbrae may represent the same sequence $1,
a_{1}, a_{2}, \ldots,$ as it happens also in probability theory for
independent and identically distributed random variables. It is mainly thanks
to these devices that the early umbral calculus has had a rigorous and simple
formal look.

At first glance, the classical umbral calculus seems just a notation for
dealing with exponential generating functions. Nevertheless, this new syntax
has given rise noteworthy computational simplifications and conceptual
clarifications in different contexts. Applications are given by Zeilberger
\cite{Zeilberger}, where generating functions are computed for many difficult
problems dealing with counting combinatorial objects. Applications to bilinear
generating functions for polynomial sequences are given by Gessel
\cite{Gessel}. Connections with wavelet theory have been investigated in
\cite{Saliani} and \cite{Shen}. In \cite{Dinardo}, the development of this
symbolic computation has produced the theory of Bell umbrae, by which the
functional composition of exponential power series has been interpreted in a
effective way. On the basis of this result, the umbral calculus has been
interpreted as a calculus of measures on Poisson algebras, generalizing
compound Poisson processes \cite{Dinardo}. A natural parallel with random
variables has been further carried out in \cite{Dinardoeurop}. In
\cite{DinardoBer}, the theory of $k$-statistics and polykays has been
completely rewritten, carrying out a unifying framework for these estimators,
both in the univariate and multivariate cases. Moreover, very fast algorithms
for computing these estimators have been carried out.

Apart from the preliminary paper of Taylor \cite{Taylor}, Sheffer sequences
have not been described in terms of umbrae. Here we complete the picture,
giving many examples and several applications.

The paper is structured as follows. Section 2 is provided for
readers unaware of the classical umbral calculus. We resume
terminology, notation and some basic definitions. In Section 3, we
introduce the notion of the adjoint of an umbra. This notion is
the key to clarify the nature of the umbral presentation of
binomial sequences. In Section 4, by introducing umbral
polynomials, we stress a feature of the classical umbral calculus,
that is the construction of new umbrae by suitable symbolic
substitutions. Section 5 is devoted to Sheffer umbrae. We
introduce the notion of Sheffer umbra from which we derive an
umbral characterization of Sheffer sequences $\{s_{n}(x)\}.$
Theorem \ref{(si)} gives the umbral version of the well-known
Sheffer identity with respect to the associated sequence. Theorem
\ref{(conrep22)} gives the umbral version of a second
characterization of Sheffer sequences, that is $\{s_{n}(x)\}$ is
said to be a Sheffer sequence with respect to a delta operator
$Q,$ when $Qs_{n}(x)=ns_{n-1}(x)$ for all $n\geq0.$ In Section 6,
the notion of Sheffer umbra is used to introducing two special
umbrae: the one associated to an umbra, whose moments are binomial
sequences, and the Appell umbra, whose moments are Appell
polynomials. We easily state their main properties by umbral
methods. In the last section, we discuss some topics to which
umbral methods can be fruitfully applied. In particular we deal
with the connection constants problem, that gives the coefficients
in expressing a sequence of polynomials $\{s_{n}(x)\}$ in terms of
a different sequence of polynomials $\{p_{n}(x)\}$ and viceversa.
We give a very simple proof of the Lagrange inversion formula by
showing that all polynomials of binomial type are represented by
Abel polynomials. This last result was proved by Rota, Shen and
Taylor in \cite{RotaTaylor}, but the authors did not make explicit
the relations among the involved umbrae, and thus have not
completely pointed out the powerfulness of the result. Moreover,
the notion of a Sheffer umbra brings to the light the umbral
connection between binomial sequences and Abel polynomials. This
allows us to give a very handy umbral expression for the Stirling
numbers of first and second type. Finally, we would like to stress
how the connection between Sheffer umbrae and Lagrange inversion
formula has smoothed the way to an umbral theory of free cumulants
\cite{free}. In closing, we provide some examples of exact
solutions of linear recursions, which benefit of an umbral
approach.
%

\section{The classical umbral calculus}

In the following, we recall terminology, notation and some basic definitions
of the classical umbral calculus, as it has been introduced by Rota and Taylor
in \cite{SIAM} and further developed in \cite{Dinardo} and \cite{Dinardoeurop}%
. An umbral calculus consists of the following data:

\begin{description}
\item[\textit{a)}] a set $A=\{\alpha,\beta, \ldots\},$ called the
\textit{alphabet}, whose elements are named \textit{umbrae};

\item[\textit{b)}] a commutative integral domain $R$ whose quotient field is
of characteristic zero;

\item[\textit{c)}] a linear functional $E,$ called \textit{evaluation},
defined on the polynomial ring $R[A]$ and taking values in $R$ such that

\begin{description}
\item[\textit{i)}] $E[1]=1;$

\item[\textit{ii)}] $E[\alpha^{i} \beta^{j} \cdots\gamma^{k}] = E[\alpha
^{i}]E[\beta^{j}] \cdots E[\gamma^{k}]$ for any set of distinct umbrae in $A$
and for $i,j,\ldots,k$ nonnegative integers (\textit{uncorrelation property});
\end{description}

\item[\textit{d)}] an element $\epsilon\in A,$ called \textit{augmentation}
\cite{Roman}, such that $E[\epsilon^{n}] = \delta_{0,n},$ for any nonnegative
integer $n,$ where
\[
\delta_{i,j} = \left\{
\begin{array}
[c]{cc}%
1, & \hbox{if $i=j$},\\
0, & \hbox{if $i \ne j$},
\end{array}
\right.  \,\, i,j \in N;
\]

\item[\textit{e)}] an element $u \in A,$ called \textit{unity} umbra
\cite{Dinardo}, such that $E[u^{n}]=1,$ for any nonnegative integer $n.$
\end{description}

A sequence $a_{0}=1,a_{1},a_{2}, \ldots$ in $R$ is umbrally represented by an
umbra $\alpha$ when
\[
E[\alpha^{i}]=a_{i}, \quad\hbox{for} \,\, i=0,1,2,\ldots.
\]
The elements $a_{i}$ are called \textit{moments} of the umbra $\alpha.$

\begin{example}
\textit{Singleton umbra.} \label{singl} \newline\textrm{The singleton umbra
$\chi$ is the umbra such that
\[
E[\chi^{1}] = 1, \quad E[\chi^{n}] = 0 \quad\hbox{for} \,\, n=2,3,\ldots.
\]
}
\end{example}

The \textit{factorial moments} of an umbra $\alpha$ are the elements
\[
a_{(n)} = \left\{
\begin{array}
[c]{ll}%
1, & n=0,\\
E[(\alpha)_{n}], & n>0,
\end{array}
\right.
\]
where $(\alpha)_{n}=\alpha(\alpha-1)\cdots(\alpha-n+1)$ is the lower factorial.

\begin{example}
\textit{Bell umbra.} \label{exbell} \newline\textrm{The Bell umbra $\beta$ is
the umbra such that
\[
E[(\beta)_{n}] = 1 \quad\hbox{for} \,\, n=0,1,2,\ldots.
\]
In \cite{Dinardo} we prove that $E[\beta^{n}]= B_{n},$ where $B_{n}$ is the
$n$-th Bell number, i.e. the number of partitions of a finite nonempty set
with $n$ elements, or the $n$-th coefficient in the Taylor series expansion of
the function $\exp(e^{t}-1).$}
\end{example}

An umbral polynomial is a polynomial $p \in R[A].$ The support of $p$ is the
set of all umbrae occurring in $p.$ If $p$ and $q$ are two umbral polynomials then

\begin{enumerate}
\item[\textit{i)}] $p$ and $q$ are \textit{uncorrelated} if and only if their
supports are disjoint;

\item[\textit{ii)}] $p$ and $q$ are \textit{umbrally equivalent} if and only
if $E[p]=E[q],$ in symbols $p\simeq q.$
\end{enumerate}


\subsection{Similar umbrae and dot-product}

The notion of similarity among umbrae comes in handy in order to manipulate
sequences such
\begin{equation}
\sum_{i=0}^{n}\left(
\begin{array}
[c]{c}%
n\\
i
\end{array}
\right)  a_{i}a_{n-i},\quad\hbox{for}\,\,n=0,1,2,\ldots\label{(eq:1)}%
\end{equation}
as moments of umbrae. The sequence (\ref{(eq:1)}) cannot be represented by
using only the umbra $\alpha$ with moments $a_{0}=1,a_{1},a_{2},\ldots.$
Indeed, $\alpha$ being correlated to itself, the product $a_{i}a_{n-i}$ cannot
be written as $E[\alpha^{i}\alpha^{n-i}].$ So, as it happens for random
variables, we need two distinct umbrae having the same sequence of moments.
Therefore, if we choose an umbra $\alpha^{\prime}$ uncorrelated with $\alpha$
but with the same sequence of moments, we have
\begin{equation}
\sum_{i=0}^{n}\left(
\begin{array}
[c]{c}%
n\\
i
\end{array}
\right)  a_{i}a_{n-i}=E\left[  \sum_{i=0}^{n}\left(
\begin{array}
[c]{c}%
n\\
i
\end{array}
\right)  \alpha^{i}(\alpha^{\prime})^{n-i}\right]
=E[(\alpha+\alpha^{\prime
})^{n}]. \label{(forment)}%
\end{equation}
Then the sequence (\ref{(eq:1)}) represents the moments of the umbra
$(\alpha+\alpha^{\prime}).$ In \cite{SIAM}, Rota and Taylor formalize this
matter by defining an equivalence relation among umbrae.

Two umbrae $\alpha$ and $\gamma$ are \textit{similar} when $\alpha^{n}$ is
umbrally equivalent to $\gamma^{n},$ for all $n=0,1,2,\ldots$ in symbols
\[
\alpha\equiv\gamma\Leftrightarrow\alpha^{n} \simeq\gamma^{n} \quad
n=0,1,2,\ldots.
\]

\begin{example}
\textit{Bernoulli umbra.} \newline\textrm{The Bernoulli umbra $\iota$ (cf.
\cite{SIAM}) satisfies the umbral equivalence $\iota+ u \equiv-\iota.$ Its
moments are the Bernoulli numbers $B_{n},$ such that
\[
\sum_{k \geq0} \left(
\begin{array}
[c]{c}%
n\\
k
\end{array}
\right)  B_{k} = B_{n}.
\]
}
\end{example}

Thanks to the notion of similar umbrae, it is possible to extend the alphabet
$A$ with the so-called \textit{auxiliary} umbrae, resulting from operations
among similar umbrae. This leads to construct a saturated umbral calculus, in
which auxiliary umbrae are handled as elements of the alphabet. It can be
shown that saturated umbral calculi exist and that every umbral calculus can
be embedded in a saturated umbral calculus \cite{SIAM}. We shall denote by the
symbol $n \boldsymbol{.} \alpha$ the \textit{dot-product} of $n$ and $\alpha,$
an auxiliary umbra (cf. \cite{SIAM}) similar to the sum $\alpha^{\prime
}+\alpha^{\prime\prime}+ \cdots+ \alpha^{\prime\prime\prime}$ where
$\{\alpha^{\prime},\alpha^{\prime\prime},\ldots,\alpha^{\prime\prime\prime}\}$
is a set of $n$ distinct umbrae, each one similar to the umbra $\alpha.$ So
the sequence in (\ref{(forment)}) is umbrally represented by the umbra $2
\boldsymbol{.} \alpha$. We assume that $0 \boldsymbol{.} \alpha$ is an
auxiliary umbra similar to the augmentation $\epsilon.$

The next statements follow from the definition of the dot-product.

\begin{proposition}
\label{(prop1)}

\begin{description}
\item[\textit{{ i)}}] If $n \boldsymbol{.} \alpha\equiv n \boldsymbol{.}
\beta$ for some integer $n \ne0,$ then $\alpha\equiv\beta;$

\item[\textit{{ii)}}] if $c \in R,$ then $n \boldsymbol{.} (c \alpha) \equiv c
(n \boldsymbol{.} \alpha) $ for any nonnegative integer $n;$

\item[\textit{{iii)}}] $n \boldsymbol{.} (m \boldsymbol{.} \alpha) \equiv(nm)
\boldsymbol{.} \alpha\equiv m \boldsymbol{.} (n \boldsymbol{.} \alpha)$ for
any two nonnegative integers $n,m;$

\item[\textit{{iv)}}] $(n+m) \boldsymbol{.} \alpha\equiv n \boldsymbol{.}
\alpha+ m \boldsymbol{.} \alpha^{\prime}$ for any two nonnegative integers
$n,m$ and any two distinct umbrae $\alpha\equiv\alpha^{\prime};$

\item[\textit{{v)}}] $(n \boldsymbol{.} \alpha+ n \boldsymbol{.} \gamma)\equiv
n \boldsymbol{.}(\alpha+\gamma)$ for any nonnegative integer $n$ and any two
distinct umbrae $\alpha$ and $\gamma.$
\end{description}
\end{proposition}

Two umbrae $\alpha$ and $\gamma$ are said to be \textit{inverse} to each other
when $\alpha+\gamma\equiv\varepsilon.$ We denote the inverse of the umbra
$\alpha$ by $-1 \boldsymbol{.} \alpha.$ Note that they are uncorrelated.
Recall that, in dealing with a saturated umbral calculus, the inverse of an
umbra is not unique, but any two umbrae inverse to any given umbra are similar.

We shall denote by the symbol $\alpha^{\boldsymbol{.} \, n}$ the
\textit{dot-power} of $\alpha$, an auxiliary umbra similar to the product
$\alpha^{\prime} \alpha^{\prime\prime}$ $\cdots\alpha^{\prime\prime\prime},$
where $\{\alpha^{\prime},\alpha^{\prime\prime},\ldots,\alpha^{\prime
\prime\prime}\}$ is a set of $n$ distinct umbrae, similar to the umbra
$\alpha.$ We assume that $\alpha^{\boldsymbol{.} \, 0}$ is an umbra similar to
the unity umbra $u.$

The next statements follow from the definition of the dot-power.

\begin{proposition}
\begin{description}

\item[\textit{{i)}}] If $c \in R,$ then $(c \alpha)^{\boldsymbol{.} \, n}
\equiv c^{ \, n} \alpha^{ \boldsymbol{.} \, n}$ for any nonnegative integer $n
\ne0;$

\item[\textit{ii)}] $(\alpha^{ \boldsymbol{.} \, n})^{\boldsymbol{.} \, m}
\equiv\alpha^{\boldsymbol{.} \, nm} \equiv(\alpha^{\boldsymbol{.} \, m})^{
\boldsymbol{.} \, n}$ for any two nonnegative integers $n,m;$

\item[\textit{iii)}] $\alpha^{\boldsymbol{.} \, (n+m)} \equiv\alpha^{
\boldsymbol{.} \, n} (\alpha^{\prime})^{\boldsymbol{.} \, m}$ for any two
nonnegative integers $n,m$ and any two distinct umbrae $\alpha\equiv
\alpha^{\prime};$

\item[\textit{iv)}] $(\alpha^{ \boldsymbol{.} \, n})^{k} \equiv(\alpha
^{k})^{\boldsymbol{.} \, n}$ for any two nonnegative integers $n,k.$
\end{description}
\end{proposition}

By the statement \textit{iv)}, the moments of $\alpha^{\boldsymbol{.} \, n}$
are:
\begin{equation}
E[(\alpha^{\boldsymbol{.} \, n})^{k}] =E[(\alpha^{k})^{\boldsymbol{.} \, n}] =
a_{k}^{n}, \quad k=0,1,2,\ldots\label{(eq:10)}%
\end{equation}
for any nonnegative integer $n.$ Hence the moments of the umbra $\alpha
^{\boldsymbol{.} \, n}$ are the $n$-th power of the moments of the umbra
$\alpha.$

Moments of $n\boldsymbol{.}\alpha$ can be expressed using the notions of
integer partitions and dot-powers. Recall that a partition of an integer $i$
is a sequence $\lambda=(\lambda_{1},\lambda_{2},\ldots,\lambda_{t}),$ where
$\lambda_{j}$ are weakly decreasing positive integers such that $\sum
_{j=1}^{t}\lambda_{j}=i.$ The integers $\lambda_{j}$ are named \textit{parts}
of $\lambda.$ The \textit{lenght} of $\lambda$ is the number of its parts and
will be indicated by $\nu_{\lambda}.$ A different notation is $\lambda
=(1^{r_{1}},2^{r_{2}},\ldots),$ where $r_{j}$ is the number of parts of
$\lambda$ equal to $j$ and $r_{1}+r_{2}+\cdots=\nu_{\lambda}.$ We use the
classical notation $\lambda\vdash i$ to denote ``$\lambda$ is a partition of
$i$''. By using an umbral version of the well-known multinomial expansion
theorem, we have
\begin{equation}
(n\boldsymbol{.}\alpha)^{i}\simeq\sum_{\lambda\vdash i}(n)_{\nu_{\lambda}%
}d_{\lambda}\alpha_{\lambda}, \label{(eq:11)}%
\end{equation}
where the sum is over all partitions $\lambda=(1^{r_{1}},2^{r_{2}},\ldots)$ of
the integer $i,$ $(n)_{\nu_{\lambda}}=0$ for $\nu_{\lambda}>n,$
\begin{equation}
d_{\lambda}=\frac{i!}{r_{1}!r_{2}!\cdots}\,\frac{1}{(1!)^{r_{1}}(2!)^{r_{2}%
}\cdots}\quad\hbox{and}\quad\alpha_{\lambda}=(\alpha_{j_{1}})^{\boldsymbol{.}%
\,r_{1}}(\alpha_{j_{2}})^{\boldsymbol{.}\,r_{2}}\cdots, \label{(eq:12)}%
\end{equation}
with $j_{i}$ distinct integers chosen in $\{1,2,\ldots,n\}.$

\subsection{The generating function of an umbra}

The formal power series
\begin{equation}
u + \sum_{n \geq1} \alpha^{n} \frac{t^{n}}{n!} \label{(gf)}%
\end{equation}
is the \textit{generating function} (g.f.) of the umbra $\alpha,$ and it is
denoted by $e^{\alpha t}.$ The notion of umbral equivalence and similarity can
be extended coefficientwise to formal power series of $R[A][[t]]$
\[
\alpha\equiv\beta\Leftrightarrow e^{\alpha t} \simeq e^{\beta t},
\]
(see \cite{Taylor1} for a formal construction). Note that any exponential
formal power series
\begin{equation}
f(t) = 1 + \sum_{n \geq1} a_{n} \frac{t^{n}}{n!} \label{(gf1)}%
\end{equation}
can be umbrally represented by a formal power series (\ref{(gf)}) in
$R[A][[t]].$ In fact, if the sequence $1,a_{1},a_{2},\ldots$ is umbrally
represented by $\alpha$ then
\[
f(t)=E[e^{\alpha t}] \quad\hbox{i.e.} \quad f(t) \simeq e^{\alpha t},
\]
assuming that we extend $E$ by linearity. We denote the formal power series in
(\ref{(gf1)}) by $f(\alpha,t)$ and we will say that $f(\alpha,t)$ is umbrally
represented by $\alpha.$ Henceforth, when no confusion occurs, we will just
say that $f(\alpha,t)$ is the g.f. of $\alpha.$ For example the g.f. of the
augmentation umbra $\epsilon$ is $f (\epsilon,t)=1$, while the g.f. of the
unity umbra $u$ is $f(u,t)=e^{t}.$ The g.f. of the singleton umbra $\chi$ is
$f(\chi,t)=1+t,$ the g.f. of the Bell umbra is $f(\beta,t)=\exp(e^{t}-1)$ and
the g.f. of the Bernoulli umbra is $f(\iota,t)=t/(e^{t}-1).$

The advantage of an umbral notation for g.f.'s is the representation of
operations among g.f.'s through symbolic operations among umbrae. For example,
the product of exponential g.f.'s is umbrally represented by a sum of the
corresponding umbrae:
\begin{equation}
f(\alpha,t)\,f(\gamma,t)\simeq e^{(\alpha+\gamma)t}\quad\hbox{with}\quad
f(\alpha,t)\simeq e^{\alpha t}\,\,\hbox{and}\,\,f(\gamma,t)\simeq e^{\gamma
t}. \label{(summation)}%
\end{equation}
Via (\ref{(summation)}), the g.f. of $n.\alpha$ is $f(\alpha,t)^{n}.$ Note
that
\begin{equation}
e^{(n\boldsymbol{.}\,\alpha)\,t}\simeq f(t)^{n}\simeq\left(  e^{\alpha
t}\right)  ^{\boldsymbol{.}\,n} \label{gfp}%
\end{equation}
Via g.f., we have \cite{Dinardo}
\begin{equation}
E[(n\boldsymbol{.}\alpha)^{i}]= \sum_{j=1}^{i}\,(n)_{j}\,B_{i,j}(a_{1}%
,a_{2},\ldots,a_{i-j+1})\quad i=1,2,\ldots\label{(gr:3bis)}%
\end{equation}
where $B_{i,j}$ are the (partial) Bell exponential polynomials (cf.
\cite{Riordan}) and $a_{i}$ are the moments of the umbra $\alpha.$

If $\alpha$ is an umbra with g.f. $f(\alpha,t),$ then
\[
e^{(-1 \boldsymbol{.} \, \alpha) \, t} \simeq\frac{1}{f(\alpha,t)}.
\]
The sum of exponential g.f.'s is umbrally represented by a disjoint sum of
umbrae. The \textit{disjoint sum} (respectively \textit{disjoint difference})
of $\alpha$ and $\gamma$ is the umbra $\eta$ (respectively $\delta$) with
moments
\[
\eta^{n} \simeq\left\{
\begin{array}
[c]{ll}%
u, & n=0\\
\alpha^{n} + \gamma^{n}, & n>0
\end{array}
\right.  \quad\left(  \hbox{respectively} \quad\delta^{n} \simeq\left\{
\begin{array}
[c]{ll}%
u, & n=0\\
\alpha^{n} - \gamma^{n}, & n>0
\end{array}
\right.  \right)  ,
\]
in symbols $\eta\equiv\alpha\dot{+} \gamma$ (respectively $\delta\equiv
\alpha\dot{-} \gamma$). By the definition, we have
\[
f(\alpha,t) \pm[f(\gamma,t) -1] \simeq e^{(\alpha\dot{\pm} \gamma)t}.
\]

\subsection{Polynomial umbrae}

The introduction of the g.f. device leads to the definition of new
auxiliary umbrae, improving the computational power of the umbral
syntax. For this purpose, we could replace $R$ by a suitable
polynomial ring having coefficients in $R$ and any desired number
of indeterminates. Then, an umbra is said to be \textit{scalar} if
the moments are elements of $R$ while it is said to be
\textit{polynomial} if the moments are polynomials. In this paper,
we deal with $R[x,y].$ In particular, we define the dot-product of
$x$ and $\alpha$ via g.f., i.e. $x \boldsymbol{.} \alpha$ is the
auxiliary umbra having g.f.
\[
e^{(x \boldsymbol{.} \alpha)} \simeq f(\alpha,t)^{x}.
\]
Proposition \ref{(prop1)} still holds, replacing $n$ with $x$ and $m$ with
$y.$

\begin{example}
\textit{Bell polynomial umbra.} \label{Bexpu} \newline\textrm{The umbra $x
\boldsymbol{.} \beta$ is the Bell polynomial umbra. Its factorial moments are
powers of $x$ and its moments are the exponential polynomials (cf.
\cite{Dinardo})
\[
(x \boldsymbol{.} \beta)_{n} \simeq x^{n} \quad\hbox{and} \quad(x
\boldsymbol{.} \beta)^{n} \simeq\sum_{k=0}^{n} S(n,k) x^{k}.
\]
Its g.f. is $f(x \boldsymbol{.} \beta, t) = \exp[x(e^{t}-1)].$}
\end{example}


\subsection{Special auxiliary umbrae}

A feature of the classical umbral calculus is the construction of new
auxiliary umbrae by suitable symbolic substitutions. In $n \boldsymbol{.}
\alpha$ replace the integer $n$ by an umbra $\gamma$. From (\ref{(gr:3bis)}),
the new auxiliary umbra $\gamma\boldsymbol{.} \alpha$ has moments
\begin{equation}
E[(\gamma\boldsymbol{.} \alpha)^{i}] = \sum_{j=1}^{i} g_{(j)} B_{i,j}
(a_{1},a_{2},\ldots, a_{i-j+1}) \quad i=1,2,\ldots\label{(gr:2bis)}%
\end{equation}
where $g_{(j)}$ are the factorial moments of the umbra $\gamma.$ The auxiliary
umbra $\gamma\boldsymbol{.} \alpha$ is called \textit{dot-product} of $\alpha$
and $\gamma.$ The g.f. $f(\gamma\boldsymbol{.} \alpha,t)$ is such that
\[
e^{(\gamma\boldsymbol{.} \, \alpha)t} \simeq[f(t)]^{\gamma} \simeq
e^{\gamma\log f(t)} \simeq g\left[  \log f(t) \right]  .
\]
Observe that $E[\gamma\boldsymbol{.} \alpha]= g_{1} \, a_{1} = E[\gamma
]\,E[\alpha.]$ The following statements hold.

\begin{proposition}
\label{prop2}

\begin{description}
\item[\textit{i)}] If $\eta\boldsymbol{.} \alpha\equiv\eta\boldsymbol{.}
\gamma$ for some umbra $\eta,$ then $\alpha\equiv\gamma;$

\item[\textit{ii)}] if $c \in R,$ then $\eta\boldsymbol{.} (c \alpha) \equiv c
(\eta\boldsymbol{.} \alpha)$ for any two distinct umbrae $\alpha$ and $\eta;$

\item[\textit{iii)}] if $\gamma\equiv\gamma^{\prime},$ then $(\alpha+\eta)
\boldsymbol{.} \gamma\equiv\alpha\boldsymbol{.} \gamma+ \eta\boldsymbol{.}
\gamma^{\prime};$

\item[\textit{iv)}] $\eta\boldsymbol{.} (\gamma\boldsymbol{.} \alpha)
\equiv(\eta\boldsymbol{.} \gamma) \boldsymbol{.} \alpha.$
\end{description}
\end{proposition}

For the proofs, see \cite{Dinardo}. Observe that from property \textit{ii)} it
follows
\begin{equation}
\alpha\boldsymbol{.} x \equiv\alpha\boldsymbol{.} (x u) \equiv x(\alpha
\boldsymbol{.} u) \equiv x\alpha. \label{(6bis)}%
\end{equation}
In the following, we recall some useful dot-products of umbrae, whose
properties have been investigated with full details in \cite{Dinardo} and
\cite{Dinardoeurop}.

\begin{example}
\textit{Exponential umbral polynomials} \label{exrp}
\newline\textrm{Suppose we replace $x$ with a generic umbra
$\alpha$ in the Bell polynomial umbra $x\boldsymbol{.}\beta.$ We
get the auxiliary umbra $\alpha\boldsymbol{.}\beta,$ whose
factorial moments are
\begin{equation}
(\alpha\boldsymbol{.}\beta)_{n}\simeq\alpha^{n}\quad n=0,1,2,\ldots.
\label{(eq:12bis)}%
\end{equation}
The moments are given by the exponential umbral polynomials (cf.
\cite{Dinardo})
\begin{equation}
(\alpha\boldsymbol{.}\beta)^{n}\simeq\Phi_{n}(\alpha)\simeq\sum_{i=0}%
^{n}S(n,i)\alpha^{i}\quad n=0,1,2,\ldots. \label{(expBell)}%
\end{equation}
The g.f. is $f(\alpha\boldsymbol{.}\beta,t)=f(e^{t}-1).$}
\end{example}

\begin{example}
\textit{$\alpha$-partition umbra} \label{exrp2} \newline\textrm{The $\alpha
$-partition umbra is the umbra $\beta\boldsymbol{.} \alpha,$ where $\beta$ is
the Bell umbra (see Example \ref{exbell}). Since the factorial moments of
$\beta$ are all equal to $1,$ equation (\ref{(gr:2bis)}) gives
\begin{equation}
E[(\beta\boldsymbol{.} \alpha)^{i}] = \sum_{j=1}^{i} B_{i,j}(a_{1}%
,a_{2},\ldots, a_{i-j+1}) = Y_{i}(a_{1},a_{2},\ldots,a_{i}) \label{(gr:4bis)}%
\end{equation}
for $i=1,2,\ldots,$ where $Y_{i}$ are the complete exponential polynomials
\cite{Riordan}. The umbra $x \boldsymbol{.} \beta\boldsymbol{.} \alpha$ is the
polynomial $\alpha$-partition umbra. Since the factorial moments of $x
\boldsymbol{.} \beta$ are powers of $x,$ equation (\ref{(gr:2bis)}) gives
\begin{equation}
E[(x \boldsymbol{.} \beta\boldsymbol{.} \alpha)^{i}] = \sum_{j=1}^{i} x^{j} \,
B_{i,j}(a_{1},a_{2},\ldots, a_{i-j+1}), \label{(gr:3ter)}%
\end{equation}
so the g.f. is $f(x \boldsymbol{.} \beta\boldsymbol{.} \alpha,t)=\exp
[x(f(t)-1)].$}
\end{example}

The $\alpha$-partition umbra $\beta\boldsymbol{.}\alpha$ plays a crucial role
in the umbral representation of the composition of exponential g.f.'s. Indeed,
the \textit{composition umbra} of $\alpha$ and $\gamma$ is the umbra
$\gamma\boldsymbol{.}\beta\boldsymbol{.}\alpha.$ The g.f. is $f(\gamma
\boldsymbol{.}\beta\boldsymbol{.}\alpha,t)=f[\gamma,f(\alpha,t)-1].$ The
moments are
\begin{equation}
E[(\gamma\boldsymbol{.}\beta\boldsymbol{.}\alpha)^{i}]=\sum_{j=1}^{i}%
g_{j}\,B_{i,j}(a_{1},a_{2},\ldots,a_{i-j+1}), \label{(momcomp)}%
\end{equation}
where $g_{j}$ and $a_{i}$ are moments of the umbra $\gamma$ and
$\alpha$ respectively. From equivalences (\ref{(eq:11)}) and
(\ref{(eq:12bis)}), we also have
\begin{equation}
(\gamma\boldsymbol{.}\beta\boldsymbol{.}\alpha)^{i}\simeq\sum_{\lambda\vdash
i}(\gamma\boldsymbol{.}\beta)_{\nu_{\lambda}}d_{\lambda}\alpha_{\lambda}%
\simeq\sum_{\lambda\vdash i}\gamma^{\nu_{\lambda}}d_{\lambda}\alpha_{\lambda}
,\label{(eq:13)}%
\end{equation}
where $d_{\lambda}$ and $\alpha_{\lambda}$ are given in
(\ref{(eq:12)}).

We denote by $\alpha^{<-1>}$ the compositional inverse of $\alpha,$ i.e. the
umbra having g.f. $f(\alpha^{<-1>},t)$ such that $f[\alpha^{<-1>}%
,f(\alpha,t)-1]=f[\alpha, f(\alpha^{<-1>},t)-1]=1+t,$ i.e. $f(\alpha
^{<-1>},t)=f^{<-1>}(\alpha,t).$ So we have
\[
\alpha\boldsymbol{.} \beta\boldsymbol{.} \alpha^{<-1>} \equiv\alpha^{<-1>}
\boldsymbol{.} \beta\boldsymbol{.} \alpha\equiv\chi.
\]
In particular for the unity umbra, we have
\begin{equation}
\beta\boldsymbol{.} u^{<-1>} \equiv u^{<-1>} \boldsymbol{.} \beta\equiv\chi,
\label{(11)}%
\end{equation}
by which the next fundamental equivalences follow
\begin{equation}
\beta\boldsymbol{.} \chi\equiv u \equiv\chi\boldsymbol{.} \beta.
\label{(fund)}%
\end{equation}
Since $\chi\boldsymbol{.} \beta\boldsymbol{.} \chi\equiv\chi\boldsymbol{.} u
\equiv\chi,$ recalling \textit{i)} of Proposition \ref{prop2}, the
compositional inverse of the singleton umbra $\chi$ is the umbra $\chi$ itself.

\begin{example}
\textit{$\alpha$-cumulant umbra} \label{(ex4)} \newline\textrm{The umbra
$\chi\boldsymbol{.} \alpha$ is the $\alpha$-cumulant umbra, having g.f.
$f(\chi\boldsymbol{.} \alpha, t)=1+\log[f(t)].$ Then, the umbra $\chi$ is the
cumulant umbra of $u,$ the umbra $u$ is the cumulant umbra of $\beta,$ the
umbra $u^{<-1>}$ is the cumulant umbra of $\chi.$ Properties of cumulant
umbrae are investigated in details in \cite{Dinardoeurop}. A special role is
held by the cumulant umbra of a polynomial Bell umbra. Indeed, as it has been
proved in \cite{Dinardoeurop}, the $(x \boldsymbol{.} \beta)$-cumulant umbra
has moments all equal to $x:$}
\begin{equation}
(\chi\boldsymbol{.} x \boldsymbol{.} \beta)^{n} \simeq x. \label{(x)}%
\end{equation}

\end{example}

\begin{example}
\textit{$\alpha$-factorial umbra} \label{(ex4)}
\newline\textrm{The umbra $\alpha\boldsymbol{.} \chi$ is the
$\alpha$-factorial umbra, since $(\alpha\boldsymbol{.} \chi)^{n}
\simeq(\alpha)_{n}$ for all nonnegative $n.$ The g.f. is
$f(\alpha.\chi,t)=f[\log(1+t)].$ By using the $\alpha$-factorial
umbra, from equivalence (\ref{(eq:11)}) we also have
\[
(\gamma\boldsymbol{.} \alpha)^{i} \simeq\sum_{\lambda\vdash i} (\gamma
)_{\nu_{\lambda}} d_{\lambda} \alpha_{\lambda} \simeq\sum_{\lambda\vdash i}
(\gamma\boldsymbol{.} \chi)^{\nu_{\lambda}} d_{\lambda} \alpha_{\lambda},
\]
where $d_{\lambda}$ and $\alpha_{\lambda}$ are given in
(\ref{(eq:12)}).}
\end{example}


\section{Adjoint umbrae}

Let $\gamma$ be an umbra with $E[\gamma] = g_{1} \ne0$ so that the g.f. $f
(\gamma,t)$ admits compositional inverse. In this section, we study some
properties of a special partition umbra, i.e the $\gamma^{<-1>}$-partition
umbra. As it will be clarified in the following, this is a key umbra in the
theory of binomial polynomials.

\begin{definition}
\label{(defadj)} The adjoint umbra of $\gamma$ is the $\gamma^{<-1>}%
$-partition umbra:
\[
\gamma^{*} = \beta\boldsymbol{.} \gamma^{<-1>}.
\]

\end{definition}

The name parallels the adjoint of an umbral operator \cite{Roman} since
$\gamma\boldsymbol{.} \alpha^{*}$ gives the umbral composition of $\gamma$ and
$\alpha^{<-1>}.$

\begin{example}
\textit{Adjoint of the singleton umbra $\chi$} \label{(adjchi)} \newline%
\textrm{\textrm{The inverse of $\chi$ is the umbra $\chi$ itself. So we have
\begin{equation}
\chi^{*} \equiv\beta\boldsymbol{.} \chi^{<-1>} \equiv\beta\boldsymbol{.}
\chi\equiv u.\label{(chiadj)}%
\end{equation}
}}
\end{example}

\begin{example}
\textit{Adjoint of the unity umbra $u$} \label{(adjchi)} \newline%
\textrm{\textrm{By virtue of equivalence (\ref{(11)}), the adjoint of the
unity umbra $u$ is
\begin{equation}
u^{*} \equiv\beta.u^{<-1>} \equiv\chi.\label{(uadj)}\\
\end{equation}
}}
\end{example}

\begin{example}
\textit{Adjoint of the Bell umbra $\beta$} \newline\textrm{\textrm{We have
$\beta^{*} \equiv u^{<-1>}.$ Indeed $\beta\boldsymbol{.} \beta\boldsymbol{.}
\beta^{<-1>} \equiv\chi$ and, taking the left-hand side dot-product by $\chi,$
we have
\[
\chi\boldsymbol{.} \beta\boldsymbol{.} \beta\boldsymbol{.} \beta^{<-1>}
\equiv\chi\boldsymbol{.} \chi\equiv u^{<-1>}.
\]
The result follows recalling equivalence (\ref{(fund)}) and $\beta
\boldsymbol{.} \beta^{<-1>} = \beta^{*}.$}}
\end{example}

The adjoint umbra has g.f.
\[
f(\gamma^{*}, t) = \exp[f^{-1}(\gamma,t)-1].
\]
In particular the adjoint of the compositional inverse of an umbra is similar
to its partition umbra, i.e.
\[
(\gamma^{<-1>})^{*} \equiv\beta\boldsymbol{.} \gamma.
\]
From the previous equivalence, we have
\begin{equation}
(\gamma^{<-1>})^{*} \boldsymbol{.} \beta\boldsymbol{.} \gamma^{<-1>}
\equiv\beta\boldsymbol{.} \gamma\boldsymbol{.} \beta\boldsymbol{.}
\gamma^{<-1>} \equiv\beta\boldsymbol{.} \chi\equiv u. \label{(1)}%
\end{equation}
and, replacing $\gamma^{<-1>}$ with $\gamma,$ we have
\begin{equation}
\gamma^{*} \boldsymbol{.} \beta\boldsymbol{.} \gamma\equiv u. \label{(2)}%
\end{equation}

\begin{example}
\textit{Adjoint of $u^{<-1>}$} \newline\textrm{\textrm{The adjoint of the
compositional inverse of the unity umbra is $(u^{<-1>})^{*} \equiv
\beta\boldsymbol{.} u \equiv\beta.$}}
\end{example}

Equivalences (\ref{(1)}) and (\ref{(2)}) may be rewritten in a more useful
way. Indeed, we have
\begin{equation}
(\gamma^{<-1>})^{*} \boldsymbol{.} \gamma^{*} \equiv\gamma^{*} \boldsymbol{.}
(\gamma^{<-1>})^{*} \equiv u. \label{(adjinv)}%
\end{equation}
Note that the dot-product of an umbra $\alpha$ with the adjoint of $\gamma$ is
the composition umbra of $\alpha$ and $\gamma^{<-1>},$ i.e. $\alpha.\gamma^{*}
\equiv\alpha\boldsymbol{.} \beta\boldsymbol{.} \gamma^{<-1>}.$ In particular
we have
\begin{equation}
\gamma\boldsymbol{.} \gamma^{*} \equiv\chi\Rightarrow\beta\boldsymbol{.}
\gamma\boldsymbol{.} \gamma^{*} \equiv\beta\boldsymbol{.} \chi\equiv u
\label{adj1}%
\end{equation}
and also
\begin{equation}
\chi\boldsymbol{.} \gamma^{*} \equiv\gamma^{<-1>} \quad\hbox{and} \quad
\chi\boldsymbol{.} (\gamma^{<-1>})^{*} \equiv\gamma. \label{adj2}%
\end{equation}

\begin{proposition}
\label{comp} The adjoint of the composition umbra of $\alpha$ and $\gamma$ is
the dot-product of the adjoints of $\gamma$ and $\alpha,$ that is
\begin{equation}
(\alpha\boldsymbol{.} \beta\boldsymbol{.} \gamma)^{*} \equiv\gamma^{*}
\boldsymbol{.} \, \alpha^{*}. \label{adjcomp}%
\end{equation}

\end{proposition}

\begin{proof}
$$(\alpha \punt \beta \punt \gamma)^* \equiv \beta \punt(\alpha \punt \beta
\punt \gamma)^{<-1>} \equiv (\beta \punt \gamma^{<-1>}) \punt
(\beta \punt \, \alpha^{<-1>}) \equiv \gamma^{*} \punt \,
\alpha^{*}.$$
\end{proof}

\section{Umbral polynomials}

Let $\{q_{n}(x)\}$ be a polynomial sequence of $R[x]$ such that $q_{n}(x)$ has
degree $n$ for any $n:$
\[
q_{n}(x)=q_{n,\,n}x^{n}+q_{n,\,n-1}x^{n-1}+\cdots+q_{n,\,0}.
\]
Moreover, let be $\alpha$ an umbra. The sequence $\{q_{n}(\alpha)\}$ consists
of umbral polynomials with support $\alpha$ such that
\[
E[q_{n}(\alpha)]=q_{n,\,n}a_{n}+q_{n,\,n-1}a_{n-1}+\cdots+q_{n,\,0}%
\]
for any nonnegative integer $n$. Now suppose $q_{0}(x)=1$ and consider an
auxiliary umbra $\eta$ such that
\[
E[\eta^{n}]=E[q_{n}(\alpha)],
\]
for any nonnegative integer $n$. In order to underline that the moments of
$\eta$ depend on those of $\alpha,$ we add the subscript $\alpha$ to the umbra
$\eta$ so that we shall write
\[
\eta_{\alpha}^{n}\simeq q_{n}(\alpha)\quad\hbox{for}\,\,n=0,1,2,\ldots.
\]
If $\alpha\equiv x.u,$ then $\eta_{x.u}$ is a polynomial umbra with moments
$q_{n}(x),$ so we shall simply denote it by $\eta_{x}.$

Let us consider some simple consequences of the notations here introduced.

\begin{proposition}
If $\eta_{x}$ is a polynomial umbra and $\alpha$ and $\gamma$ are umbrae both
scalar either polynomial, then
\[
\eta_{\alpha}\equiv\eta_{\gamma}\Leftrightarrow\alpha\equiv\gamma.
\]

\end{proposition}

\begin{proof}
For any nonnegative integer $n,$ there exist constants $c_{n, \,
k}, k=0,1,\ldots,n$ such that $x^n=\sum_{k=0}^n c_{n,k} \,
q_k(x).$ Since $\eta_\alpha \equiv \eta_\gamma,$ then $q_k(\alpha)
\simeq q_k(\gamma)$ for all nonnegative integers $k$ and so for
all nonnegative integers $n$ we have
$$\alpha^n \simeq \sum_{k=0}^n c_{n, \,k} \, q_k(\alpha) \simeq \sum_{k=0}^n
c_{n, \, k} \, q_k(\gamma) \simeq \gamma^n.$$  The other direction
of the proof is straight forward.
\end{proof}

\begin{proposition}
\label{(unic)} If $\eta_{x}$ and $\zeta_{x}$ are polynomial umbrae and
$\alpha$ is an umbra either scalar or polynomial, then
\[
\eta_{\alpha}\equiv\zeta_{\alpha}\Leftrightarrow\eta_{x}\equiv\zeta_{x}.
\]

\end{proposition}

\begin{proof}
Suppose $\eta_\alpha \equiv \zeta_\alpha.$ Let $\{q_n(x)\}$ be the
moments of $\eta_x$ and let $\{z_n(x)\}$ be the moments of
$\zeta_x.$ For all nonnegative integers $n$, there exist constants
$c_{n, \, k}, k=0,1,\ldots,n$ such that $q_n(x)=\sum_{k=0}^n c_{n,
\, k} \, z_k(x).$ Being $q_n(\alpha) \simeq z_n(\alpha)$ for all
nonnegative integers $n$, we have $c_{n, \, k}=\delta_{n, \, k}$
for $k=0,1,\ldots,n$ by which we have $q_n(x)= z_n(x).$ The other
direction of the proof is straight forward.
\end{proof}

\begin{proposition}
If $\eta_{x}$ is a polynomial umbra, $\{\alpha_{i}\}_{i=1}^{n}$ are
uncorrelated scalar umbrae and $\{w_{i}\}_{i=1}^{n}$ are some weights in $R,$
then
\begin{equation}
\eta_{\dot{+}_{i=1}^{n} \chi.w_{i}.\beta.\alpha_{i}} \equiv\dot{+}_{i=1}^{n}
\chi.w_{i}.\beta.\eta_{\alpha_{i}}. \label{(poldis)}%
\end{equation}

\end{proposition}

\begin{proof}
Suppose $E[\eta_x^m]=q_{m}(x)=\sum_{k=0}^m q_{m, \, k} x^k.$ Then
for all nonnegative integers $m$,  we have
\begin{eqnarray*}
q_m( \dot{+}_{i=1}^n \chi.w_i.\beta.\alpha_i) & \simeq &
\sum_{k=0}^m q_{m,k} ( \dot{+}_{i=1}^n \chi.w_i.\beta.\alpha_i)^k
\simeq  \sum_{k=0}^m q_{m,k} \left( \sum_{i=1}^n w_i \alpha_i^k\right)
\\
& \simeq & \sum_{i=1}^n w_i \left(\sum_{k=0}^m q_{m,k}
\alpha_i^k\right) \simeq   \dot{+}_{i=1}^n
\chi.w_i.\beta.q_m(\alpha_i),
\end{eqnarray*}
due to equivalence (\ref{(x)}). The result follows by observing
that
$$\dot{+}_{i=1}^n \chi.w_i.\beta.q_m(\alpha_i) \simeq \left[\dot{+}_{i=1}^n
\chi.w_i.\beta.\eta_{\alpha_i})\right]^m.$$
\end{proof}
We achieve the proof of the next corollary choosing $w_{i}=1$ for
$i=1,2,...,n$ in equivalence (\ref{(poldis)}).

\begin{corollary}
If $\eta_{x}$ is a polynomial umbra and $\{\alpha_{i}\}_{i=1}^{n}$ are
uncorrelated scalar umbrae, then
\[
\eta_{\dot{+}_{i=1}^{n} \alpha_{i}} \equiv\dot{+}_{i=1}^{n} \eta_{\alpha_{i}%
}.
\]

\end{corollary}


\section{Sheffer sequences}

In this section we give the definition of Sheffer umbra by which we recover
fundamental properties of Sheffer sequences. In the following let us $\alpha$
and $\gamma$ be scalar umbrae, with $g_{1}=E[\gamma] \ne0.$

\begin{definition}
\label{(defshef)} A polynomial umbra $\sigma_x$ is said to be a
Sheffer umbra for $(\alpha,\gamma)$ if
\begin{equation}
\sigma_x \equiv(-1 \boldsymbol{.} \alpha+ x \boldsymbol{.} u)
\boldsymbol{.}
\gamma^{*}, \label{(conrep)}%
\end{equation}
where $\gamma^{*}$ is the adjoint umbra of $\gamma.$
\end{definition}

In the following, we denote a Sheffer umbra by $\sigma_{x}^{(\alpha,\gamma)}$
in order to make explicit the dependence on $\alpha$ and $\gamma.$

We note that if $\alpha$ has g.f. $f(\alpha, t)$ and $\gamma$ has g.f.
$f(\gamma, t),$ the g.f of $(-1 \boldsymbol{.} \alpha+ x \boldsymbol{.} u)
\boldsymbol{.} \gamma^{*}$ is the composition of $f(\gamma^{<-1>},
t)=f^{-1}(\gamma,t)$ and $f(- 1 \boldsymbol{.} \alpha+ x \boldsymbol{.} u, t)
= e^{xt}/f(\alpha, t),$ i.e.
\begin{equation}
f( \sigma^{(\alpha,\gamma)}_{x}, t) =
\frac{1}{f[\alpha,f^{-1}(\gamma,t)-1]}
e^{x \, [f^{-1}(\gamma,t)-1]}. \label{gf}%
\end{equation}

\begin{theorem}
[The Expansion Theorem]\label{expthm} If $\sigma^{(\alpha,\gamma)}_{x}$ is a
Sheffer umbra for $(\alpha,\gamma),$ then
\begin{equation}
\eta\equiv\alpha+ \sigma^{(\alpha,\gamma)}_{\eta} \boldsymbol{.}
\beta\boldsymbol{.} \gamma\label{(expansionth)}%
\end{equation}
for any umbra $\eta.$
\end{theorem}

\begin{proof}
Replacing $x$ with $\eta$ in equivalence (\ref{(conrep)}),  we
obtain
$$\sigma^{(\alpha,\gamma)} _{\eta} \equiv (- 1 \punt \alpha + \eta) \punt
\gamma^{*}.$$ Take the right dot product with $\beta \punt \gamma$
of both sides, then
$$\sigma^{(\alpha,\gamma)} _{\eta} \punt \, \beta \punt \gamma \equiv (- 1 \punt \alpha + \eta) \punt
\gamma^{*} \punt \, \beta \punt \gamma \equiv (- 1 \punt \alpha +
\eta).$$ The result follows adding $\alpha$ to both sides of the
previous equivalence.
\end{proof}

\begin{theorem}
\label{(th0)} Let $\{s_{n}(x)\}$ be the moments of a Sheffer umbra
$\sigma^{(\alpha,\gamma)}_{x}.$ The polynomial sequence $\{s_{n}(x)\}$ is the
unique polynomial sequence such that:
\begin{equation}
s_{n}(\alpha+ k \boldsymbol{.} \gamma) \simeq(k \boldsymbol{.} \chi)^{n}
\quad\forall n,k \geq0. \label{(sheffer)}%
\end{equation}

\end{theorem}

\begin{proof}
From equivalence (\ref{(conrep)}), when $x$ is replaced by $\alpha
+ k \punt \gamma,$ we have
\begin{equation}
(-1 \punt \alpha+\alpha+k \punt \gamma) \punt \gamma^* \equiv k
\punt \gamma \punt \gamma^* \equiv k \punt \chi \quad \forall \, k
\geq 0 \label{(sheffer1)}
\end{equation}
which gives (\ref{(sheffer)}). The uniqueness follows from
Proposition \ref{(unic)}.
\end{proof}We explicitly note that any Sheffer umbra is uniquely determined by
its moments evaluated at $0,$ since via equivalence
(\ref{(conrep)}) we have
\[
\sigma_{0}^{(\alpha,\gamma)}\equiv-1\boldsymbol{.}\alpha\boldsymbol{.}%
\gamma^{\ast}.
\]
In Theorem \ref{(si)} we will prove that the moments of the umbra $\sigma
_{x}^{(\alpha,\gamma)}$ satisfy the Sheffer identity.

\begin{example}
\textit{Power polynomials.} \textrm{\textrm{Choosing as umbra
$\alpha$ the umbra $\epsilon$ and as umbra $\gamma$ the umbra
$\chi,$ from equivalence (\ref{(conrep)}) we have
\[
\sigma^{(\epsilon,\chi)}_{x} \equiv x \boldsymbol{.} u,
\]
since $\chi^{*} \equiv u.$ So the sequence of polynomials
$\{x^{n}\}$ is a Sheffer sequence, being moments of the Sheffer
umbra $\sigma^{(\epsilon ,\chi)}_{x}.$}}
\end{example}

\begin{example}
\label{C.P} \textit{Poisson-Charlier polynomials.} \textrm{\textrm{The
Poisson-Charlier polynomials are
\[
c_{n}(x;a)=\frac{1}{a^{n}}\sum_{k=0}^{n}\left(
\begin{array}
[c]{c}%
n\\
k
\end{array}
\right)  (-a)^{n-k}(x)_{k},
\]
hence
\[
c_{n}(x;a)\simeq\left\{  \frac{x\boldsymbol{.}\chi-a}{a}\right\}  ^{n},
\]
recalling that $(x \punt \chi)^{k}\simeq(x)_{k}.$ Denoting by
$\omega_{x,a}$ the polynomial umbra whose moments are
$c_{n}(x;a),$ we have
\[
\omega_{x,a}\equiv\frac{x \punt \chi-a}{a}.
\]
The umbra $\omega_{x,a}$ is called the Poisson-Charlier polynomial umbra. We
show that $\omega_{x,a}$ is a Sheffer umbra for $(a\boldsymbol{.}\beta
,\chi\boldsymbol{.}a\boldsymbol{.}\beta).$ Indeed
\[
\omega_{x,a}\equiv\frac{x\boldsymbol{.}\chi-a\boldsymbol{.}\beta
\boldsymbol{.}\chi}{a}\equiv\lbrack-1\boldsymbol{.}(a\boldsymbol{.}%
\beta)+x\boldsymbol{.}u]\boldsymbol{.}\frac{\chi}{a}.
\]
So the Poisson-Charlier polynomial umbra $\omega_{x,a}$ is a Sheffer umbra,
being
\[
\frac{\chi}{a}\equiv(\chi\boldsymbol{.}a\boldsymbol{.}\beta)^{\ast}.
\]
}}
\end{example}

\bigskip

\begin{theorem}
[The Sheffer identity]\label{(si)} A polynomial umbra $\sigma_{x}$
is a Sheffer umbra if and only if there exists an umbra $\eta,$
provided with a compositional inverse, such that
\begin{equation}
\sigma_{x+y}\equiv\sigma_{x}+y\boldsymbol{.}\eta^{\ast}.\label{(shefid)}%
\end{equation}

\end{theorem}

\begin{proof}
Let $\sigma_x$ be a Sheffer umbra for $(\alpha, \gamma).$ By
Definition \ref{(defshef)}, we have
\begin{eqnarray*}
\sigma^{(\alpha,\gamma)}_{x+y}  & \equiv &  (-1 \punt \alpha + x
\punt u + y
\punt u) \punt \gamma^* \\
& \equiv & (-1 \punt \alpha + x \punt u) \punt \gamma^* + y \punt
\gamma^* \equiv \sigma^{(\alpha,\gamma)}_x + y \punt \gamma^*.
\end{eqnarray*}
Viceversa, set  $x=0$ in equivalence (\ref{(shefid)}). We have
\begin{equation}
\sigma_y \equiv \sigma_0 + y \punt \eta^*. \label{(conrep1)}
\end{equation}
On the other hand, known the moments of $\eta,$  there exists an
umbra $\alpha$ such that $-1 \punt \alpha \punt \eta^* \equiv
\sigma_0.$ Then the polynomial umbra $\sigma_y$ is a Sheffer
umbra, being $\sigma_y \equiv (-1 \punt \alpha+y \punt u) \punt
\eta^*.$
\end{proof}Equivalence (\ref{(shefid)}) gives the well-known Sheffer identity,
because by using the binomial expansion we have
\begin{equation}
s_{n}(x+y)=\sum_{k=0}^{n}\left(
\begin{array}
[c]{c}%
n\\
k\\
\end{array}
\right)  \,s_{k}(x)\,p_{n-k}(y)\label{(shefid1)}%
\end{equation}
where $s_{n}(x+y)=E[\sigma_{x+y}^{n}],$
$s_{k}(x)=E[\sigma_{x}^{k}]$ and
$p_{n-k}(y)=E[(y\boldsymbol{.}\eta^{\ast})^{n-k}].$ In the next
section, we will prove that the moments of umbrae such
$y\boldsymbol{.}\eta^{\ast}$ are binomial sequences.

\begin{corollary}
If $\sigma^{(\alpha,\gamma)}_{x}$ is a Sheffer umbra for $(\alpha,\gamma),$
then
\[
\sigma^{(\alpha,\gamma)}_{\eta+\zeta}
\equiv\sigma^{(\alpha,\gamma)}_{\eta} + \zeta \boldsymbol{.}
\gamma^{*},
\]
where $\eta$ and $\zeta$ are umbrae.
\end{corollary}

\begin{theorem}
\label{(conrep22)} A polynomial umbra $\sigma_{x}$ is a Sheffer umbra if and
only if there exists an umbra $\eta,$ provided with compositional inverse,
such that
\begin{equation}
\sigma_{\eta+x \boldsymbol{.} u} \equiv\chi+ \sigma_{x}. \label{(conrep2)}%
\end{equation}

\end{theorem}

\begin{proof}
If $\sigma_x$ is a Sheffer umbra for $(\alpha,\gamma),$ then
$$\sigma^{(\alpha,\gamma)}_{\gamma+x \punt u}  \equiv (-1 \punt \alpha + x \punt u +
\gamma) \punt \gamma^* \equiv (-1 \punt \alpha + x \punt u) \punt
\gamma^* + \chi,$$ recalling that $\gamma \punt \gamma^* \equiv
\chi.$ Then equivalence (\ref {(conrep2)}) follows from
equivalence (\ref{(conrep)}) choosing as umbra $\eta$ the umbra
$\gamma.$ Viceversa, let $\sigma_x$ be a polynomial umbra such
that equivalence (\ref{(conrep2)}) holds for some umbra $\eta,$
with $E[\eta] = g_1 \ne 0.$ Set $x=0$ in equivalence
(\ref{(conrep2)}). We have
$$\sigma_{\eta} \equiv \chi + \sigma_0 \Rightarrow
\sigma_{\eta} \punt \beta \punt \eta \equiv \chi \punt \beta \punt
\eta + \sigma_0 \punt \beta \punt \eta \equiv \eta + \sigma_0
\punt \beta \punt \eta.$$ Known the moments of $\eta,$ there
exists an umbra $\alpha$ such that $-1 \punt \alpha \punt \eta^*
\equiv \sigma_0$ so that
$$\sigma_{\eta} \punt \beta \punt \eta \equiv \eta  -1 \punt \alpha.$$
Due to equivalence (\ref{(expansionth)}), also  the Sheffer umbra
for $(\alpha,\eta)$ is such that $\sigma^{(\alpha,\eta)}_{\eta}
\punt \beta \punt \eta \equiv \eta  -1 \punt \alpha,$ therefore
$\sigma_{\eta} \equiv \sigma^{(\alpha,\eta)}_{\eta}$ and
$\sigma_{x} \equiv \sigma^{(\alpha,\eta)}_{x}$ by Proposition
\ref{(unic)}.
\end{proof}

\begin{corollary}
\label{(deriv)} A polynomial umbra $\sigma_{x}$ is a Sheffer umbra if and only
if there exists an umbra $\eta,$ provided with compositional inverse, such
that
\[
\sigma^{k}_{\eta+x.u} \simeq\sigma^{k}_{x} + k \, \sigma^{k-1}_{x}
\quad\hbox{\rm
for} \,\, k=1,2,\ldots.
\]

\end{corollary}

\begin{proof}
Take the $k$-th moment of both sides in equivalence
(\ref{(conrep2)}).
\end{proof}
Now suppose $s_{n}(x)=\sum_{k=0}^{n} s_{n,k} \, x^{k}$ be the moments of a
Sheffer umbra for $(\alpha,\gamma)$ and $r_{n}(x)$ be the moments of a Sheffer
umbra for $(\eta,\zeta).$ The umbra
\[
[-1 \boldsymbol{.} \alpha+(-1 \boldsymbol{.} \eta+x \boldsymbol{.} u)
\boldsymbol{.} \zeta^{*}] \boldsymbol{.} \gamma^{*} \equiv[-1 \boldsymbol{.}
(\alpha\boldsymbol{.} \beta\boldsymbol{.} \zeta+\eta)+x \boldsymbol{.} u]
\boldsymbol{.} \zeta^{*} \boldsymbol{.} \gamma^{*}%
\]
has moments umbrally equivalent to
\begin{equation}
\sum_{k=0}^{n} s_{n,k} r_{k}(x), \label{(shpocp)}%
\end{equation}
i.e. the Roman-Rota umbral composition $s_{n}(\mathbf{r}(x)).$ So we have
proved the following theorem.

\begin{theorem}
[Umbral composition and Sheffer umbrae]The polynomials given in
(\ref{(shpocp)}) are moments of the Sheffer umbra for $(\alpha
\boldsymbol{.} \beta \boldsymbol{.} \zeta+\eta,\gamma
\boldsymbol{.} \beta \boldsymbol{.} \zeta).$
\end{theorem}

Two Sheffer sequences are said to be inverse of each other if and only if
their Roman-Rota umbral composition (\ref{(shpocp)}) gives the sequence
$\{x^{n}\}.$

\begin{corollary}
[Inverse of Sheffer sequences]\label{invsh} The sequence of moments
corresponding to the Sheffer umbra for $(-1 \boldsymbol{.} \alpha
\boldsymbol{.} \gamma^{*},\gamma^{<-1>})$ are inverses of the sequence of
moments corresponding to the Sheffer umbra for $(\alpha,\gamma).$
\end{corollary}

\begin{proof}
Set
$$\zeta \equiv \gamma^{<-1>} \quad \hbox{and} \quad \eta \equiv -1 \punt \alpha
\punt \gamma^{*},$$ in $(\alpha \punt \beta \punt
\zeta+\eta,\gamma \punt \beta \punt \zeta).$ The corresponding
Sheffer umbra has moments umbrally equivalent to $x^n.$
\end{proof}

\section{Two special Sheffer umbrae}

In this section, we study two special classes of Sheffer umbrae: the
associated umbra and the Appell umbra. The associated umbrae are polynomial
umbrae whose moments $\{p_{n}(x)\}$ satisfies the well-known binomial
identity
\begin{equation}
p_{n}(x+y) = \sum_{k=0}^{n} \left(
\begin{array}
[c]{c}%
n\\
k\\
\end{array}
\right)  \, p_{k}(x) \, p_{n-k}(y) \label{(binom1)}%
\end{equation}
for all $n = 0,1,2,\ldots.$ Every sequence of binomial type is a Sheffer
sequence but most Sheffer sequences are not of binomial type. The concept of
binomial type has applications in combinatorics, probability, statistics, and
a variety of other fields. The Appell umbrae are polynomial umbrae whose
moments $\{p_{n}(x)\}$ satisfies the identity
\begin{equation}
\frac{d}{dx} p_{n}(x) = n p_{n-1}(x) \quad n = 1,2,\ldots. \label{(Appell)}%
\end{equation}
Among the most notable Appell sequences, besides the trivial example
$\{x^{n}\},$ are the Hermite polynomials, the Bernoulli polynomials, and the
Euler polynomials.

\subsection{Associated umbrae}

Let us consider a Sheffer umbra for the umbrae $(\epsilon,\gamma),$ where
$\gamma$ has compositional inverse and $\epsilon$ is the augmentation umbra.

\begin{definition}
\label{(th1bis)} A polynomial umbra $\sigma_x$ is said to be the
associated umbra of $\gamma$ if
\[
\sigma_x \equiv x \boldsymbol{.} \gamma^{*},
\]
where $\gamma^{*}$ is the adjoint umbra of $\gamma.$
\end{definition}

The g.f. of $x \boldsymbol{.} \gamma^{*}$ is
\[
f(x \boldsymbol{.} \gamma^{*},t)= e^{x [f^{-1}(\gamma,t)-1]},
\]
because in equation (\ref{gf}) we have $f(\alpha,t) =
f(\epsilon,t) = 1.$ The expansion theorem for associated umbrae
(cfr. Theorem \ref{expthm}) is
\begin{equation}
\eta\equiv\eta\boldsymbol{.} \gamma^{*} \boldsymbol{.} \beta\boldsymbol{.}
\gamma. \label{(expass)}%
\end{equation}

\begin{theorem}
\label{(212)} An umbra $\sigma^{(\alpha,\gamma)}_{x}$ is a Sheffer umbra for
$(\alpha,\gamma)$ if and only if $\sigma^{(\alpha,\gamma)}_{\alpha+x
\boldsymbol{.} u}$ is the umbra associated to $\gamma.$
\end{theorem}

\begin{proof}
If $\sigma^{(\alpha,\gamma)}_x$ is a Sheffer umbra for
$(\alpha,\gamma)$ then
$$\sigma^{(\alpha,\gamma)}_{\alpha+x \punt u}  \equiv
(-1 \punt \alpha + \alpha \punt u + x \punt u) \punt \gamma^*,$$
by which
$$\sigma^{(\alpha,\gamma)}_{\alpha+ x \punt u} \equiv x \punt \gamma^*.$$
From Definition \ref{(th1bis)}, the umbra
$\sigma^{(\alpha,\gamma)}_{\alpha+x \punt u}$ is the umbra
associated to $\gamma.$ Viceversa, let $\eta_x$ be a polynomial
umbra such that
$$\eta_{\alpha + x \punt u} \equiv x \punt \gamma^*.$$
Replacing $x$ with $k \punt \gamma,$ we have
$$\eta_{\alpha + k \punt \gamma} \equiv k \punt \gamma \punt \gamma^* \equiv k
\punt \chi.$$ The result follows by equivalence
(\ref{(sheffer1)}).
\end{proof}
We will say that a polynomial sequence $\{p_{n}(x)\}$ is \textit{associated}
to an umbra $\gamma$ if and only if
\[
p_{n}(x) \simeq(x \boldsymbol{.} \gamma^{*})^{n}, \quad n=0,1,2,\ldots
\]
or
\begin{equation}
p_{n}(k \punt \gamma) \simeq(k \boldsymbol{.} \chi)^{n} \quad
n,k=0,1,2,\ldots.
\label{(associated)}%
\end{equation}

\begin{theorem}
[Umbral characterization of associated sequences]\label{ucp} The sequence
$\{p_{n}(x)\}$ is associated to the umbra $\gamma$ if and only if:
\begin{align}
&  p_{n}(\epsilon) \simeq\epsilon^{n} \quad\hbox{for} \,\, n = 0,1,2,\ldots
\label{(ass1)}\\
&  p_{n}(\gamma+x \boldsymbol{.} u ) \simeq p_{n}(x) + n \, p_{n-1}(x)
\quad\hbox{for} \,\, n = 1,2,\ldots. \label{(ass2)}%
\end{align}

\end{theorem}

\begin{proof}
If the sequence $\{p_n(x)\}$ is associated to $\gamma$ then
$$p_n(\epsilon) \simeq (\epsilon \punt \gamma^*)^n \simeq \epsilon^n$$
for all $n = 0,1,2,\ldots.$ Equivalence (\ref{(ass2)}) follows
from Corollary \ref{(deriv)} choosing as umbra $\alpha$ the umbra
$\epsilon.$ Viceversa, if equivalences (\ref{(ass1)}) and
(\ref{(ass2)}) hold, we prove by induction that the sequence
$\{p_n(x)\}$ satisfies (\ref{(associated)}). Indeed, by
equivalence (\ref{(ass1)}) we have
$$p_n(0 \punt \gamma) \simeq p_n(\epsilon) \simeq \left\{ \begin{array}{ll}
0, & \,\,\, n = 1,2,\ldots,j, \\
1, & \,\,\, n = 0.
\end{array} \right.$$
Suppose that equivalence  (\ref{(associated)}) holds for $k=m$
\begin{equation}
p_n(m \punt \gamma) \simeq (m \punt \chi)^n \quad n=0,1,2,\ldots.
\label{(hi)}
\end{equation}
By equivalence (\ref{(ass2)}), we have
$$p_n[(m+1) \punt \gamma) ] \simeq p_n(\gamma + m \punt \gamma) \simeq p_n(m
\punt \gamma) + n p_{n-1}(m \punt \gamma) \quad n=1,2,\ldots.$$
Due to induction hypothesis (\ref{(hi)}), we have
$$p_n[(m+1) \punt \gamma] \simeq  (m \punt \chi)^n + n (m \punt \chi)^{n-1} \simeq
(\chi + m \punt \chi)^n \simeq [(m+1)\punt \chi]^n \quad
n=1,2,\ldots.$$ Since the sequence $\{p_n(x)\}$ verifies
(\ref{(associated)}), it is associated to $\gamma.$
\end{proof}

\begin{theorem}
[The binomial identity]The sequence $\{p_{n}(x)\}$ is associated to the umbra
$\gamma$ if and only if
\begin{equation}
p_{n}(x+y) = \sum_{k=0}^{n} \left(
\begin{array}
[c]{c}%
n\\
k\\
\end{array}
\right)  p_{k}(x) \, p_{n-k}(y) \label{(binom)}%
\end{equation}
for all $n = 0,1,2,\ldots.$
\end{theorem}

\begin{proof}
If the sequence $\{p_n(x)\}$ is associated to the umbra $\gamma,$
then identity (\ref{(binom)}) follows from the property
$$(x+y) \punt \gamma^* \equiv x \punt \gamma^* + y \punt \gamma^*.$$
Viceversa, suppose the sequence $\{p_n(x)\}$ satisfies identity
(\ref{(binom)}).  Let $\eta_x$ be a polynomial umbra such that
$$E[\eta^n_x] = p_n(x).$$
By identity  (\ref{(binom)}), we have
\begin{equation}
\eta_{x+y} \equiv \eta_x + \eta^{\prime}_y \label{(binom1)}
\end{equation}
with $\eta_x$ similar to $\eta^{\prime}_x$ and uncorrelated. In
particular, if we replace $y$ with $\epsilon$ in equivalence (\ref
{(binom1)}), then $\eta_x \equiv \eta_x +
\eta^{\prime}_{\epsilon}$  and hence $\eta^{\prime}_{\epsilon}
\equiv \epsilon.$ So the polynomials $\{p_n(x)\}$ are such that
$p_n(\epsilon) \simeq \epsilon^n,$ i.e. they satisfy equivalence
(\ref {(ass1)}). By induction on equivalence (\ref{(binom1)}), we
have
$$\eta_{\scriptscriptstyle{{\underbrace{x + \cdots + x}_k}}} \equiv \underbrace{\eta_x + \cdots +
\eta^{\prime}_x}_k$$ where the polynomial umbrae on the right-hand
side are uncorrelated and similar to $\eta_x.$ If the $x$'s are
replaced by uncorrelated umbrae similar to any umbra $\gamma,$
provided of compositional inverse, then
$$\eta_{k \punt \gamma} \equiv k \punt \eta_{\gamma}.$$
Since $E[\gamma] \ne 0,$ we can choose an umbra $\gamma$ such that
$$\eta_{\gamma} \equiv \chi$$
thus
\begin{equation}
\eta_{k \punt \gamma} \equiv k \punt \chi \label{(ass3)}
\end{equation}
and the result follows from equivalences (\ref{(ass3)}) and (\ref
{(associated)}).
\end{proof}

\begin{example}
\textrm{\textrm{The umbra $x \boldsymbol{.} u$ is associated to the umbra
$\chi.$ Indeed, the polynomial sequence $\{x^{n} \}$ is associated to the
adjoint umbra $\chi^{*} \equiv\beta\boldsymbol{.} \chi^{<-1>} \equiv u.$ The
g.f. is
\[
f(x \boldsymbol{.} u,t) = \sum_{k \geq0} \frac{x^{k}}{k!} t^{k}%
\]
and the binomial identity becomes
\[
(x+y)^{n} = \sum_{k=0}^{n} \left(
\begin{array}
[c]{c}%
n\\
k\\
\end{array}
\right)  x^{k} \, y^{n-k}.
\]
}}
\end{example}

\begin{example}
\textrm{\textrm{The umbra $x \boldsymbol{.} u^{*} \equiv x
\boldsymbol{.} \chi$ is associated to the umbra $u.$ The
associated polynomial sequence is given by $\{ (x \punt \chi)^{n}
\} \simeq\{ (x)_{n} \},$ see Example \ref{(ex4)}. The g.f. is
\[
f(x \boldsymbol{.} u^{*}, t) = (1+t)^{x}%
\]
and the binomial identity becomes
\[
(x+y)_{n} = \sum_{k=0}^{n} \left(
\begin{array}
[c]{c}%
n\\
k\\
\end{array}
\right)  (x)_{k} \, (y)_{n-k}.
\]
}}
\end{example}

\begin{example}
\textrm{\textrm{The umbra $x \boldsymbol{.} (u^{<-1>})^{*} \equiv
x \boldsymbol{.} \beta$ is associated to the umbra $u^{<-1>}.$ The
associated polynomial sequence is given by $\{ (x \boldsymbol{.}
\beta) ^{n} \} \simeq\{\Phi_{n}(x)\},$ where $\{\Phi_{n}(x)\}$ are
the exponential polynomials (\ref{(expBell)}).}}
\end{example}

Now, suppose $\{p_{n}(x)\}$ be the polynomial sequence associated to an umbra
$\gamma$ with g.f. $f(\gamma,t)$ and $\{q_{n}(x)\}$ be the polynomial sequence
associated to an umbra $\zeta$ with g.f. $f(\zeta, t)$, i.e.
\[
p_{n}(x) \simeq(x \boldsymbol{.} \beta\boldsymbol{.} \gamma^{<-1>})^{n}
\quad\hbox{\rm and} \quad q_{n}(x) \simeq(x \boldsymbol{.} \beta\boldsymbol{.}
\zeta^{<-1>})^{n}%
\]
for $n=0,1,2,\ldots.$ On the other hand, due to Proposition \ref{comp}, we
have
\[
(x \boldsymbol{.} \beta\boldsymbol{.} \gamma^{<-1>}) \boldsymbol{.}
\beta\boldsymbol{.} \zeta^{<-1>} \equiv x \boldsymbol{.} (\beta\boldsymbol{.}
\gamma^{<-1>}) \boldsymbol{.} (\beta\boldsymbol{.} \zeta^{<-1>}) \equiv x
\boldsymbol{.} \gamma^{*} \boldsymbol{.} \zeta^{*} \equiv x \boldsymbol{.}
(\zeta\boldsymbol{.} \beta\boldsymbol{.} \gamma)^{*}.
\]
Following Roman-Rota notation, the umbra $x \boldsymbol{.} (\zeta
\boldsymbol{.} \beta\boldsymbol{.} \gamma)^{*}$ has moments
\begin{equation}
q_{n}(\mathbf{p}(x))=\sum_{k=0}^{n} q_{n,k} p_{k}(x). \label{(compass)}%
\end{equation}
This proves the next theorem.

\begin{theorem}
[Umbral composition of associated sequences]The polynomial sequence
(\ref{(compass)}) is associated to the compositional umbra of $\zeta$ and
$\gamma.$
\end{theorem}

\begin{corollary}
[Inverse of associated sequences]\label{invass} The polynomial sequence
associated to the umbra $\gamma^{<-1>}$ is inverse of the polynomial sequence
associated to the umbra $\gamma.$
\end{corollary}

\begin{proof}
Choosing as umbra $\gamma$ the umbra $\zeta^{<-1>}$ in the
previous theorem, we have
$$x \punt \zeta^* \punt \beta \punt \zeta \equiv x \punt u.$$
\end{proof}

\begin{example}
\textrm{\textrm{Choose as umbra $\zeta$ the umbra $u.$ Then $x \boldsymbol{.}
u^{*} \boldsymbol{.} \beta\boldsymbol{.} u \equiv x \boldsymbol{.}
\chi\boldsymbol{.} \beta\equiv x\boldsymbol{.} u$ and so
\[
x^{n} = \sum_{k=0}^{n} S(n,k) (x)_{k},
\]
which is the well-known formula giving powers in terms of lower factorials.}}
\end{example}

Finally, via Proposition 3.1, it is easy to prove the following recurrence
formula
\[
(x \boldsymbol{.} \gamma^{*})^{n+1} \simeq x \, \gamma^{<-1>} \, [(x+\chi)
\boldsymbol{.} \gamma^{*}]^{n}.
\]

\subsection{Appell umbrae}

In this section, we consider a second kind of special Sheffer umbra, that is a
Sheffer umbra for $(\alpha,\chi).$

\begin{definition}
\label{(th2bis)} A polynomial umbra $\sigma_x$ is said to be the
Appell umbra of $\alpha$ if
\[
\sigma_x \equiv-1 \boldsymbol{.} \alpha+ x \boldsymbol{.} u.
\]

\end{definition}

By equivalence (\ref{gf}), being $f(\chi,t)=1+t,$ the g.f. of $(-1
\boldsymbol{.} \alpha+ x \boldsymbol{.} u)$ is
\[
f(-1 \boldsymbol{.} \alpha+ x \boldsymbol{.} u, t) = \frac{1}{f(\alpha,t)} \,
e^{x t}.
\]
We will say that a polynomial sequence $\{p_{n}(x)\}$ is an Appell sequence if
and only if
\[
p_{n}(x) \simeq(-1 \boldsymbol{.} \alpha+ x \boldsymbol{.} u)^{n}, \quad
n=0,1,2,\ldots.
\]
The expansion theorem for Appell umbrae (cfr. Theorem \ref{expthm}) easily
follows by observing that
\[
\eta\equiv\alpha+ (-1 \boldsymbol{.} \alpha+ \eta).
\]

\begin{theorem}
[The Appell identity]The polynomial umbra $\sigma_{x}$ is an Appell umbra for
some umbra $\alpha$ if and only if
\[
\sigma_{x+y} \equiv\sigma_{x} + y.
\]

\end{theorem}

\begin{proof} The result follows immediately, choosing as umbra $\gamma$ the
singleton umbra $\chi$ in the Sheffer identity.
\end{proof}

\begin{corollary}
\label{(deriveq1)} A polynomial umbra $\sigma_{x}$ is an Appell umbra for some
umbra $\alpha$ if and only if
\begin{equation}
\sigma^{n}_{\chi+x \boldsymbol{.} u} \simeq\sigma^{n}_{x} + n \sigma^{n-1}(x).
\label{(deriveq)}%
\end{equation}

\end{corollary}

\begin{proof}
The result follows from Corollary \ref{(deriv)}, choosing as umbra
$\gamma$ the singleton umbra $\chi$
\end{proof}
Roughly speaking, Corollary \ref{(deriveq1)} says that, when in the Appell
umbra we replace $x$ by $\chi+ x \boldsymbol{.} u,$ the umbra $\chi$ acts as a
derivative operator. Indeed by the binomial expansion, we have
\begin{align*}
(-1 \boldsymbol{.} \alpha+ \chi+ x \boldsymbol{.} u)^{n}  &  \simeq\sum_{k
\geq0} \left(
\begin{array}
[c]{c}%
n\\
k
\end{array}
\right)  (- 1 \boldsymbol{.} \alpha)^{n-k} \, (\chi+ x \boldsymbol{.} u)^{k}\\
&  \simeq\sum_{k \geq0} \left(
\begin{array}
[c]{c}%
n\\
k
\end{array}
\right)  (- 1 \boldsymbol{.} \alpha)^{n-k} \, [(x \boldsymbol{.} u)^{k} + k (x
\boldsymbol{.} u) ^{k-1}]\\
&  \simeq(-1 \boldsymbol{.} \alpha+ x \boldsymbol{.} u)^{n} + \sum_{k \geq0}
\left(
\begin{array}
[c]{c}%
n\\
k
\end{array}
\right)  (- 1 \boldsymbol{.} \alpha)^{n-k} \, D_{x}[(x \boldsymbol{.} u)^{k}].
\end{align*}
So equivalence (\ref{(deriveq)}) umbrally expresses equation (\ref{(Appell)}).

\begin{theorem}
[The Multiplication Theorem]For any constant $c \in R,$ we have
\[
-1 \boldsymbol{.} \alpha+ (c \, x) \boldsymbol{.} u \equiv c(-1 \boldsymbol{.}
\alpha+ x \boldsymbol{.} u) + (c- 1) \boldsymbol{.} \alpha.
\]

\end{theorem}

\begin{proof} We have
$$c(-1 \punt \alpha + x \punt u) + (c-1) \punt \alpha \equiv -c \punt \alpha +
(c \, x) \punt u + c \punt \alpha -1 \punt \alpha \equiv -1 \punt
\alpha+ (c \, x) \punt u.$$
\end{proof}

\begin{example}
[Bernoulli polynomials]\textrm{\textrm{The Appell umbra for the inverse of the
Bernoulli umbra is $\iota+ x \boldsymbol{.} u.$ From the binomial expansion,
its moments are the Bernoulli polynomials
\[
E[(\iota+ x \boldsymbol{.} u)^{n}] = \sum_{k \geq0} \left(
\begin{array}
[c]{c}%
n\\
k
\end{array}
\right)  B_{n-k} \, x^{k},
\]
where $B_{n}$ are the Bernoulli numbers.}}
\end{example}


\section{Topics in umbral calculus}

In this section, we apply the language of umbrae to some topics which benefit
from this approach. In particular we discuss the well-known connection
constants problem, the Lagrange inversion formula, and we solve some
recurrence relations to give an indication of the effectiveness of the
method.

\subsection{The connection constants problem}

The connection constants problem consists in determining the connection
constants $c_{n,k}$ in the expression
\[
s_{n}(x)=\sum_{k=0}^{n} c_{n,k} r_{k}(x),
\]
where $s_{n}(x)$ and $r_{n}(x)$ are sequences of polynomials. When $s_{n}(x)$
and $r_{n}(x)$ are Sheffer sequences, umbrae provide an easy solution to this
problem. Indeed, suppose $\eta_{x}$ be a polynomial umbra such that
\[
E[\eta^{n}_{x}]=q_{n}(x)=\sum_{k=0}^{n} c_{n,k} x^{k}.
\]
The theorem we are going to prove states that $\eta_{x}$ is a Sheffer umbra
whenever $s_{n}(x)$ and $r_{n}(x)$ are Sheffer sequences.

\begin{theorem}
If $\eta_{x}$ is a polynomial umbra such that
\begin{equation}
(-1 \boldsymbol{.} \alpha+ x \boldsymbol{.} u) \boldsymbol{.} \gamma^{*}
\equiv\eta_{(-1 \boldsymbol{.} \delta+x \boldsymbol{.} u) \boldsymbol{.}
\zeta^{*}}, \label{(cc1)}%
\end{equation}
then $\eta_{x}$ is a Sheffer umbra for $((\delta-1\boldsymbol{.} \alpha)
\boldsymbol{.} \zeta^{*},\gamma\boldsymbol{.} \beta\boldsymbol{.} \zeta
^{<-1>}).$
\end{theorem}

Note that equivalence (\ref{(cc1)}) is the way to transform the Sheffer umbra
$(-1 \boldsymbol{.} \delta+ x \boldsymbol{.} u) \boldsymbol{.} \zeta^{*}$ in
the Sheffer umbra $(-1 \boldsymbol{.} \alpha+ x \boldsymbol{.} u)
\boldsymbol{.} \gamma^{*}$ by using the polynomial umbra $\eta_{x}.$
\begin{proof}
In equivalence (\ref{(cc1)}), replace $x$ with $\delta+x \punt
\beta \punt \zeta.$ The right-hand side of equivalence
(\ref{(cc1)}) becomes
$$\eta_{[-1 \punt \delta+(\delta+x \punt \beta \punt \zeta) \punt u] \punt
\zeta^*} \equiv \eta_{x \punt \beta \punt \zeta \punt\zeta^*}
\equiv \eta_{x \punt u}$$ due to equivalence (\ref{adj1}), whereas
the left-hand side gives
\begin{equation}
[-1 \punt \alpha+(\delta+x \punt \beta \punt \zeta) \punt u] \punt
\gamma^* \equiv  (\delta-1 \punt \alpha) \punt \gamma^* + x \punt
\beta \punt \zeta \punt \gamma^*. \label{(app)}
\end{equation}
By Proposition \ref{comp}, we have
$$x  \punt \beta \punt \zeta \punt \gamma^* \equiv x \punt (\zeta^{<-1>})^*
\punt \gamma^* \equiv x \punt (\gamma \punt \beta \punt
\zeta^{<-1>})^*$$ and by equivalence (\ref{(adjinv)})
$$(\delta-1 \punt \alpha) \punt \gamma^* \equiv (\delta-1 \punt \alpha) \punt u
\punt \gamma^* \equiv (\delta-1 \punt \alpha) \punt \zeta^* \punt
(\zeta^{<-1>})^* \punt \gamma^*.$$ Replacing $(\delta-1 \punt
\alpha) \punt \gamma^*$ and $x  \punt \beta \punt \zeta \punt
\gamma^*$ in equivalence (\ref{(app)}), we have
\begin{equation}
[-1 \punt \alpha+(\delta+x \punt \beta \punt \zeta) \punt u]\punt
\gamma^* \equiv [(\delta-1 \punt \alpha) \punt \zeta^* + x \punt
u] \punt (\gamma \punt \beta \punt \zeta^{<-1>})^*
\label{(cambiobase1)}
\end{equation}
and so equivalence (\ref{(cc1)}) returns
\begin{equation}
\eta_x \equiv [(\delta-1 \punt \alpha) \punt \zeta^* + x \punt
u]\punt(\gamma \punt \beta \punt \zeta^{<-1>})^*.
\label{(cambiobase)}
\end{equation}
The result follows from Definition \ref{(defshef)}.
\end{proof}
To get an explicit expression of the connection constants, we can use
equivalence (\ref{(eq:13)}) to expand the $n$-th moment of $\eta_{x}$ in
(\ref{(cambiobase)})
\[
\eta^{n}_{x} \simeq\sum_{\lambda\vdash n} [(\delta-1 \punt \alpha)
\punt \zeta^{*} + x \punt u]^{\nu_{\lambda}} \, d_{\lambda} \,
(\zeta \punt \beta \punt \gamma^{<-1>})_{\lambda}.
\]
The connection constants $c_{n,k}$ are the coefficient of $x^{k}$ in the
previous equivalence.

\begin{example}
\textrm{\textrm{Consider
\[
c_{n}(x;b) =\sum_{k=0}^{n} c_{n,k} \, c_{k}(x;a)
\]
where $c_{n}(x;a)$ and $c_{n}(x;b)$ are Poisson-Charlier polynomials. As shown
in Example \ref{C.P}, $c_{n}(x;b)$ are the moments of a Sheffer umbra for $(b
\boldsymbol{.} \beta,\chi\boldsymbol{.} b \boldsymbol{.} \beta)$ and
$c_{n}(x;a)$ are the moments of a Sheffer umbra for $(a \boldsymbol{.}
\beta,\chi\boldsymbol{.} a \boldsymbol{.} \beta).$ By equivalence
(\ref{(cambiobase)}), we have
\[
\eta_{x} \equiv[(a \boldsymbol{.} \beta-b \boldsymbol{.} \beta) \boldsymbol{.}
(\chi\boldsymbol{.} a \boldsymbol{.} \beta) ^{*} + x \boldsymbol{.} u]
\boldsymbol{.} [(\chi\boldsymbol{.} b \boldsymbol{.} \beta) \boldsymbol{.}
\beta\boldsymbol{.} (\chi\boldsymbol{.} a \boldsymbol{.} \beta)^{<-1>}]^{*}%
\]
or, equivalently, via equivalence (\ref{(cambiobase1)})
\[
\eta_{x} \equiv(-b \boldsymbol{.} \beta+ a \boldsymbol{.} \beta+ x
\boldsymbol{.} \beta\boldsymbol{.} \chi\boldsymbol{.} a \boldsymbol{.} \beta)
\boldsymbol{.} (\chi\boldsymbol{.} b \boldsymbol{.} \beta)^{*}.
\]
Observe that $x \boldsymbol{.} \beta\boldsymbol{.} \chi\boldsymbol{.} a
\boldsymbol{.} \beta\equiv xa \boldsymbol{.} \beta,$ so
\begin{equation}
(-b \boldsymbol{.} \beta+ a \boldsymbol{.} \beta+ x \boldsymbol{.}
\beta\boldsymbol{.} \chi\boldsymbol{.} a \boldsymbol{.} \beta) \boldsymbol{.}
(\chi\boldsymbol{.} b \boldsymbol{.} \beta)^{*} \equiv(a-b+xa) \boldsymbol{.}
\beta\boldsymbol{.} (\chi\boldsymbol{.} b \boldsymbol{.} \beta)^{*}.\label{CP}%
\end{equation}
Moreover $(\chi\boldsymbol{.} b \boldsymbol{.} \beta)^{*} \equiv
(\chi\boldsymbol{.} b \boldsymbol{.} \beta\boldsymbol{.} u) ^{*} \equiv u^{*}
\boldsymbol{.} (b\chi)^{*} \equiv\chi\boldsymbol{.} (b\chi)^{*},$ from which
\[
\beta\boldsymbol{.} (\chi\boldsymbol{.} b \boldsymbol{.} \beta)^{*}
\equiv\beta\boldsymbol{.} \chi\boldsymbol{.} (b\chi) ^{*} \equiv(b\chi)^{*}.
\]
Since $\chi^{*} \equiv u$ then $(b\chi)^{*} \equiv u/b$ and
\[
\beta\boldsymbol{.} (\chi\boldsymbol{.} b \boldsymbol{.} \beta)^{*} \equiv
u/b.
\]
Replacing this last result in equivalence (\ref{CP}), we have
\[
\eta_{x} \equiv\frac{a-b+xa}{b}.
\]
The binomial expansion gives
\[
c_{n,k}= \left(
\begin{array}
[c]{c}%
n\\
k
\end{array}
\right)  \left(  \frac{a}{b} \right)  ^{n} \left(  1 - \frac{b}{a} \right)
^{n-k}.
\]
}}
\end{example}


\subsection{Abel polynomials and Lagrange inversion formula}

Abel polynomials play a leading role in the theory of associated
sequences of polynomials. The main result of this section is the
proof that any sequence of binomial type can be represented as
Abel polynomials, heart of the paper \cite{RotaTaylor}. The proof
given in \cite{RotaTaylor} was a hybrid, based both on the early
Roman-Rota version of the umbral calculus and the last version,
introduced by Rota-Taylor. Here, we give a very simple proof by
introducing the notion of the derivative of an umbra.

\begin{definition}
\label{(defder)} The derivative umbra $\alpha_{\scriptscriptstyle D}$ of an
umbra $\alpha$ is the umbra whose moments are:
\[
(\alpha_{\scriptscriptstyle D})^{n} \simeq\partial_{\alpha} \alpha^{n} \simeq
n \alpha^{n-1} \quad\hbox{\rm for} \,\, n=1,2,\ldots.
\]

\end{definition}

We have
\[
f(\alpha_{\scriptscriptstyle D},t)=1 +t \, f(\alpha,t),
\]
since
\[
e^{\alpha_{\scriptscriptstyle D} t} \simeq u + \sum_{n \ge1} n \alpha^{n-1}
\frac{t^{n}}{n!} \simeq u + t e^{\alpha t}.
\]
In particular, we have
\begin{equation}
\left(  e^{\alpha_{\scriptscriptstyle D} t} - u \right)  ^{\boldsymbol{.} k}
\simeq t^{k} \left(  e^{\alpha t} \right)  ^{\boldsymbol{.} k} \simeq t^{k}
e^{(k \boldsymbol{.} \alpha) t}. \label{(genfun1)}%
\end{equation}
Note that $E[\alpha_{\scriptscriptstyle D}]=1$.

\begin{example}
[Singleton umbra]\label{sin} \textrm{\textrm{The singleton umbra $\chi$ is the
derivative umbra of the augumentation umbra $\epsilon,$ that is $\epsilon
_{\scriptscriptstyle D} \equiv\chi.$}}
\end{example}

\begin{example}
[Bernoulli umbra]\label{bu} \textrm{\textrm{We have $u \equiv(-1
\boldsymbol{.} \iota)_{\scriptscriptstyle D}.$ Indeed, we have
\[
f (-1 \boldsymbol{.} \iota, t) = \frac{e^{t} - 1}{t}%
\]
so that
\[
f [(-1 \boldsymbol{.} \iota)_{\scriptscriptstyle D}, t] = 1 + t \frac{e^{t} -
1}{t} = e^{t} = f(u,t).
\]
}}
\end{example}

\begin{example}
[Bernoulli-factorial umbra]\label{bu1} \textrm{\textrm{We have $u^{<-1>}
\equiv(\iota\boldsymbol{.} \chi)_{\scriptscriptstyle
D}.$ Indeed
\[
f (\iota\boldsymbol{.} \chi, t) = \frac{\log(1+t)}{e^{\log(1+t)}-1} =
\frac{\log(1+t)} {t}%
\]
so that
\[
f [(\iota\boldsymbol{.} \chi)_{\scriptscriptstyle D}, t] =1 + t \frac
{\log(1+t)}{t} = 1 + \log(1+t) = f(u^{<-1>},t).
\]
}}
\end{example}

\begin{theorem}
[Abel representation of binomial sequences]\label{Abel} If
$\gamma$ is an umbra provided  with a compositional inverse, then
\begin{equation}
(x \punt \gamma_{\scriptscriptstyle D}^{\ast})^{n}\simeq x(x-n
\punt \gamma)^{n-1},\quad
n=1,2,\dots \label{(Abel)}%
\end{equation}
for all $x\in R.$
\end{theorem}
In the following, we refer to polynomials $x(x-n \punt
\gamma)^{n-1}$ as \textit{umbral Abel polynomials}.
\begin{proof} On the basis of Theorem \ref{ucp}, the result follows showing that
umbral Abel polynomials are associated to the umbra $\gamma,$ i.e.
showing that such polynomials satisfy equivalences (\ref{(ass1)})
and (\ref{(ass2)}). Since  $\epsilon (\epsilon -n \punt
\gamma)^{n-1} \simeq \epsilon^n,$ equivalences (\ref{(ass1)}) are
satisfied. Moreover, it is easy to check by simple calculations
that
$$(x \punt u +\gamma_{\scriptscriptstyle D})^n - x^n \simeq n (x \punt u
+\gamma)^{n-1}, \quad n=1,2,\ldots$$ and more in general
$$p(x \punt u +\gamma_{\scriptscriptstyle D}) - p(x) \simeq \frac{d}{d x} p(x
\punt u +\gamma)$$ for any polynomial $p(x) \in R[A][x].$ In
particular for $p_n(x)=x(x-n \punt \gamma)^ {n-1},$ we have
$$p_n(x \punt u +\gamma_{\scriptscriptstyle D}) -p_n(x) \simeq \frac{d}{d x} p_n
(x\punt u+\gamma).$$ We state equivalences (\ref{(ass2)}) by
proving that
$$\frac{d}{d x} p_n(x \punt u+\gamma) \simeq n p_{n-1}(x).$$
To this aim, we have
$$\frac{d}{d x} p_n(x) \simeq n (x-n \punt \gamma)^{n-2}(x-1 \punt \gamma)$$
so
$$\frac{d}{d x} p_n(x \punt u +\gamma) \simeq n x (x-(n-1)\punt \gamma)^{n-2} \simeq
n p_{n-1}(x).$$
\end{proof}According to Corollary \ref{invass}, the inverses of umbral Abel
polynomials with respect to the Roman-Rota umbral composition are
the moments of $x \punt \beta \punt \gamma_{\scriptscriptstyle
D}.$ An umbral expression of the inverses of umbral Abel
polynomials will be given in Corollary \ref{ccc}.

Since the g.f. of $x \punt \gamma_{\scriptscriptstyle D}^{\ast}$
is $\exp [x\,(f^{<-1>}(\gamma_{\scriptscriptstyle D},t)-1)],$ we
have
\[
\exp[x\,(f^{<-1>}(\gamma_{\scriptscriptstyle D},t)-1)]\simeq1+\sum_{k\geq
1}\frac{t^{k}}{k!}[x(x-k\boldsymbol{.}\gamma)]^{k-1},
\]
which is the g.f. of umbral Abel polynomials.

Theorem \ref{Abel} includes the well-known Transfer Formula. In the following
we state various results usually derived by Transfer Formula. We start with
the Lagrange inversion formula.

\begin{corollary}
For any umbra $\gamma,$ we have
\begin{equation}
(\gamma_{\scriptscriptstyle D}^{<-1>})^{n} \simeq(-n \boldsymbol{.}
\gamma)^{n-1}, \quad n=1,2,\ldots. \label{(InvL)}%
\end{equation}

\end{corollary}

\begin{proof}
Since $\chi \punt \beta \equiv u$ then $\gamma_{\scriptscriptstyle
D}^{<-1>} \equiv \chi \punt \beta \punt \gamma_{\scriptscriptstyle
D}^{<-1>}.$ From equivalence (\ref{(Abel)}), with $x$ replaced by
$\chi,$ we have
$$(\gamma_{\scriptscriptstyle D}^{<-1>})^n \simeq  (\chi \punt \beta \punt
\gamma_{\scriptscriptstyle D}^{<-1>})^n \simeq \chi(\chi-n \punt
\gamma)^{n-1}, \quad n=1,2, \ldots.$$ Being $\chi^{k+1} \simeq 0$
for $k=1,2,\ldots,n-1,$ we have
$$\chi(\chi-n \punt \gamma)^{n-1} \simeq \sum_{k=0}^{n-1} \left( \begin{array}
{c}
n-1 \\
k \\
\end{array} \right) \chi^{k+1} (n \punt \gamma)^{n-1-k} \simeq (n \punt \gamma)^
{n-1}$$ by which equivalence (\ref{(InvL)}) follows.
\end{proof}

\begin{remark}
\label{gLIF} \textrm{\textrm{Theorem \ref{Abel} is referred to normalized
binomial polynomials, i.e. sequences $\{p_{n}(x)\}$ such that $p_{1}(x)$ is
monic. This is why the Lagrange inversion formula (\ref{(InvL)}) refers to
umbrae having first moment equal to $1$. More in general, if one would
consider umbrae having first moment different from zero, one step more is
necessary. By way of an example, we do this for the Lagrange inversion
formula.}}

\textrm{For any umbra $\gamma$ such that $E[\gamma]=g_{1} \ne0$, there exists
\footnote{In this case, in the setting of the umbral calculus $R$ must be a
field.} an umbra $\alpha$ such that $\gamma/ g_{1} \equiv\alpha
_{\scriptscriptstyle D}.$ Indeed such an umbra $\alpha$ has moments
\[
\alpha^{n-1} \simeq\frac{\gamma^{n}}{n \, g_{1}^{n}} \quad n=1,2,\ldots
\]
and g.f. $f(\alpha,t) = [f(\gamma,t/g_{1})-1]/t.$ In particular $g_{1}
\alpha\equiv\overline{\gamma},$ where $\overline{\gamma}$ is the umbra
introduced in \cite{Dinardo} having moments
\[
E[\overline{\gamma}^{\, n}] = \frac{g_{n+1}}{g_{1} (n + 1)} \quad
n=0,1,\ldots.
\]
with $E[\gamma^{n}]=g_{n}, n=1,2,\ldots.$ Multiplying for $g_{1}^{n-1}$ both
sides of equivalence (\ref{(InvL)}), written for the umbra $\alpha,$ we have
\[
g_{1}^{n-1}(\alpha_{\scriptscriptstyle D}^{<-1>})^{n} \simeq g_{1}^{n-1} (-n
\boldsymbol{.} \alpha)^{n-1} \simeq[-n \boldsymbol{.} (g_{1} \alpha)]^{n-1}
\simeq(-n \boldsymbol{.} \overline{\gamma} \,)^{n-1}%
\]
for $n=1,2,\ldots.$ Being $\gamma\boldsymbol{.} \beta\boldsymbol{.}
\gamma^{<-1>} \equiv\chi,$ recalling equivalence (\ref{(6bis)}) and
$\beta\boldsymbol{.} \chi\equiv u,$ we have
\[
\frac{\gamma}{g_{1}} g_{1} \boldsymbol{.} \beta\boldsymbol{.} \gamma^{<-1>}
\Leftrightarrow\frac{\gamma}{g_{1}} \boldsymbol{.} g_{1} \boldsymbol{.}
\beta\boldsymbol{.} \gamma^{<-1>} \equiv\chi\Leftrightarrow\frac{\gamma}%
{g_{1}} \boldsymbol{.} \beta\boldsymbol{.} \chi\boldsymbol{.} g_{1}
\boldsymbol{.} \beta\boldsymbol{.} \gamma^{<-1>} \equiv\chi.
\]
So, we have
\[
\left(  \frac{\gamma}{g_{1}} \right)  ^{<-1>} \equiv\chi\boldsymbol{.} g_{1}
\boldsymbol{.} \beta\boldsymbol{.} \gamma^{<-1>}%
\]
and from equivalence (\ref{(eq:13)}) we have
\[
\left[  \left(  \frac{\gamma}{g_{1}} \right)  ^{<-1>}\right]  ^{n} \simeq
(\chi\boldsymbol{.} g_{1} \boldsymbol{.} \beta\boldsymbol{.} \gamma
^{<-1>})^{n} \simeq g_{1} (\gamma^{<-1>})^{n}.
\]
Finally, being
\[
g_{1}^{n-1}(\alpha_{\scriptscriptstyle D}^{<-1>})^{n} \simeq g_{1}^{n-1}
\left[  \left(  \frac{\gamma}{g_{1}} \right)  ^{<-1>}\right]  ^{n} \simeq
g_{1}^{n} (\gamma^{<-1>})^{n}%
\]
we have
\begin{equation}
\gamma^{\boldsymbol{.} n} (\gamma^{<-1>})^{n} \simeq(-n \boldsymbol{.}
\overline{\gamma} \,)^{n-1} \label{(InvL2)}%
\end{equation}
This last equivalence is the generalized Lagrange inversion formula.}
\end{remark}

The Lagrange inversion formula (\ref{(InvL)}) allows us to express Stirling
numbers of first kind in terms of the Bernoulli umbra. An analogous result was
proved by Rota and Taylor in \cite{SIAM} via N\"{o}rlund sequences.

\begin{proposition}
If $\iota$ is the Bernoulli umbra and $s(n,1)=(-1)^{n-1}(n- 1)!$
$n=1,2,\ldots$ are the Stirling numbers of first kind, then
\begin{equation}
s(n,1) \simeq(n \punt \iota)^{n-1}, \quad n=1,2,\ldots. \label{(inv:3)}%
\end{equation}

\end{proposition}

\begin{proof}
From Example \ref{bu}, we have $u^{<-1>} \equiv (-1 \punt
\iota)_{\scriptscriptstyle D}^{<-1>},$ where $\iota$ is the
Bernoulli umbra. From equivalence (\ref{(InvL)}), we have
$$(u^{<-1>})^n \simeq \left[(-1 \punt \iota)_{\scriptscriptstyle D}^{<-1>)} \right]^n
\simeq [-n \punt (-1 \punt \iota)]^{n-1} \simeq (n \punt
\iota)^{n-1}, n=1,2,\ldots.$$ Equivalence (\ref{(inv:3)}) follows
recalling that $1+\log(1+t)$ is the g.f. of $u^{<-1>}$ and that
$$1+\log(1+t)=1+\sum_{n=1}^{\infty} s(n,1) \frac{t^n}{n!}$$
with $s(k,1)$ the Stirling numbers of first kind.
\end{proof}
One more application of Theorem \ref{Abel} is the proof of the following
proposition, giving a property of Abel polynomials, known as \textit{Abel
identity}.

\begin{proposition}
[Abel identity]If $\gamma\in A,$ then
\begin{equation}
(x+y)^{n} \simeq\sum_{k \geq0} \left(
\begin{array}
[c]{c}%
n\\
k
\end{array}
\right)  y(y-k \boldsymbol{.} \gamma)^{k-1} (x + k \boldsymbol{.}
\gamma)^{n-k}. \label{(polexp22)}%
\end{equation}

\end{proposition}

\begin{proof}
Recall that
$$e^{(y \punt \beta \punt \gamma_{\scriptscriptstyle D}) t} \simeq \sum_{k \geq
0} y^k \frac{(e^{\gamma_{\scriptscriptstyle D} t}-u)^{\punt
k}}{k!}.$$ Replace $y$ by $y \punt \gamma_{\scriptscriptstyle
D}^*.$ Since $\gamma_ {\scriptscriptstyle D}^* \punt \beta \punt
\gamma_{\scriptscriptstyle D} \equiv u$ we have
\begin{eqnarray}
e^{y t} & \simeq & \sum_{k \geq 0} (y \punt
\gamma_{\scriptscriptstyle D}^*)^k
\frac{(e^{\gamma_{\scriptscriptstyle D} t}-u)^{\punt k}}{k!}
\simeq  \sum_{k \geq 0} (y \punt \gamma_{\scriptscriptstyle D}^*)^k \frac{t^k
e^{(k \punt \gamma)t}}{k!} \label{(expansionth22)} \\
& \simeq & \sum_{k \geq 0} y(y-k \punt \gamma)^{k-1} \frac{t^k
e^{(k \punt \gamma) t}}{k!} \label{(expansionth3)}
\end{eqnarray}
where the second equivalence in (\ref{(expansionth22)}) follows
from (\ref {(genfun1)}), and equivalence (\ref{(expansionth3)})
follows from Theorem \ref{Abel}. Multiplying both sides by $e^{x
t},$ we have
\begin{equation}
e^{(x+y) t} \simeq \sum_{k \geq 0} y(y-k \punt \gamma)^{k-1}
\frac{t^k e^{(x+k \punt \gamma)t}}{k!}. \label{(exp1)}
\end{equation}
Since $(x+y)^n = D^{(n)}_t[e^{(x+y) t}]_{t=0},$ where
$D^{(n)}_t[\cdot]_{t=0}$ is the $n$-th derivative with respect to
$t$ evaluated at $t=0,$ the result follows taking the $n$-th
derivative with respect to $t$ of the  right-hand side of
(\ref{(exp1)}) and evaluating it at $t=0.$ Indeed, by using the
binomial property of the derivative operator we have
$$D^{(n)}_t[t^k \, e^{(x+k \punt \gamma)t}] \simeq \sum_{j=0}^k \left( \begin
{array}{c}
n \\
j \end{array} \right) D^{(j)}_t[t^k] \, D_t^{(n-j)}[e^{(x+k \punt
\gamma) t}]$$ and, setting $t=0,$ we have
$$D^{(n)}_t[t^k e^{(x + k \punt \gamma)t}]_{t=0} \simeq (n)_k (x + k \punt
\gamma)^{n-k}$$ by which equivalence (\ref{(polexp22)}) follows.
\end{proof}Setting $x=0$ in equivalence (\ref{(polexp22)}), we obtain
\begin{equation}
y^{n}\simeq\sum_{k\geq0}\left(
\begin{array}
[c]{c}%
n\\
k
\end{array}
\right)  (k\boldsymbol{.}\gamma)^{n-k}y(y-k\boldsymbol{.}\gamma)^{k-1}.
\label{(exp22)}%
\end{equation}
This proves the following polynomial expansion theorem in terms of Abel polynomials.

\begin{proposition}
\label{pe} If $p(x) \in R[x],$ then
\[
p(x) \simeq\sum_{k \geq0}p^{(k)}(k \boldsymbol{.} \gamma) \frac{y(y-k
\boldsymbol{.} \gamma)^{k-1}}{k!}.
\]

\end{proposition}

The following corollary gives the umbral expression of the Bell exponential
polynomials in (\ref{(gr:3ter)}).

\begin{corollary}
[Umbral representation of Bell exponential polynomials]\label{ccc} For all
nonnegative $n,$ we have
\begin{equation}
(x \boldsymbol{.} \beta\boldsymbol{.} \gamma_{\scriptscriptstyle D})^{n}
\simeq\sum_{k \geq0} \left(
\begin{array}
[c]{c}%
n\\
k\\
\end{array}
\right)  (k \boldsymbol{.} \gamma)^{n-k} x^{k}. \label{(exp23)}%
\end{equation}

\end{corollary}

\begin{proof}
From equivalence (\ref{(exp22)}) and  by using Theorem
\ref{Abel}, we have
\begin{equation}
y^n \simeq \sum_{k \geq 0} \left( \begin{array}{c}
n \\
k \end{array} \right) (k \punt \gamma)^{n-k} (y \punt \gamma_
{\scriptscriptstyle D}^*)^{k}. \label{(exp2)}
\end{equation}
Replace $y$ with $x \punt \beta \punt
\gamma_{\scriptscriptstyle_D}.$ We have
$$
(x \punt \beta \punt \gamma_{\scriptscriptstyle D})^n \simeq
\sum_{k \geq 0} \left( \begin{array}{c}
n \\
k \end{array} \right) (k \punt \gamma)^{n-k} (x \punt \beta \punt
\gamma_ {\scriptscriptstyle_D} \punt
\gamma_{\scriptscriptstyle_D}^*)^{k} \simeq \sum_{k \geq 0} \left(
\begin{array}{c}
n \\
k \end{array} \right) (k \punt \gamma)^{n-k} x^{k},
$$
by which the result follows immediately recalling equivalence
(\ref{(adjinv)}) for $ \gamma_{\scriptscriptstyle_D}.$
\end{proof}
The generalization of equivalence (\ref{(exp23)}) to umbrae $\gamma$ with
first moment $g_{1}$ different from zero can be stated by using the same
arguments given in Remark \ref{gLIF}:
\begin{align*}
(x \boldsymbol{.} \beta\boldsymbol{.} \gamma)^{n}  &  \simeq\gamma
^{\boldsymbol{.} n} (x \boldsymbol{.} \beta\boldsymbol{.} \alpha
_{\scriptscriptstyle D})^{n} \simeq\sum_{k \geq0} \left(
\begin{array}
[c]{c}%
n\\
k\\
\end{array}
\right)  \gamma^{\boldsymbol{.} n} (k \boldsymbol{.} \alpha)^{n-k} x^{k}\\
&  \simeq\sum_{k \geq0} \left(
\begin{array}
[c]{c}%
n\\
k\\
\end{array}
\right)  \gamma^{\boldsymbol{.} k} [k \boldsymbol{.} (g_{1} \alpha)]^{n-k}
x^{k} \simeq\sum_{k \geq0} \left(
\begin{array}
[c]{c}%
n\\
k\\
\end{array}
\right)  \gamma^{\boldsymbol{.} k} [k \boldsymbol{.} \overline{\gamma
}\,]^{n-k} x^{k}.%
\end{align*}

\begin{example}
\label{StII}\textbf{(Stirling numbers of second kind)} \newline%
\textrm{\textrm{In equivalence (\ref{(exp23)}), choose $-1
\boldsymbol{.} \iota$ as umbra $\gamma,$ then $x \boldsymbol{.}
\beta\boldsymbol{.} (-1 \boldsymbol{.}
\iota)_{\scriptscriptstyle_{D}} \equiv x \boldsymbol{.} \beta$
(see Example \ref{bu}). Comparing equivalence (\ref{(exp23)}) with
$(x \punt \beta) ^{n} \simeq\sum_{k \geq0} S(n,k) x^{k},$ where
$\{S(n,k)\}$ are Stirling numbers of second kind (see
\cite{Dinardo}), we have
\[
S(n,k) \simeq\left(
\begin{array}
[c]{c}%
n\\
k
\end{array}
\right)  (-k \punt \iota)^{n-k} \quad k=0,1,\ldots,n.
\]
This last equivalence was already proved by Rota and Taylor in \cite{SIAM}
through a different approach.}}
\end{example}

The umbral version of Stirling numbers of first kind is given in the following proposition.

\begin{proposition}
If $\{s(n,k)\}$ are Stirling numbers of first kind, then
\begin{equation}
s(n,k) \simeq\left(
\begin{array}
[c]{c}%
n\\
k
\end{array}
\right)  (k \punt \iota \punt \chi)^{n-k} \quad k=0,1,\ldots,n. \label{(eq:35)}%
\end{equation}

\end{proposition}

\begin{proof}
Recalling Example \ref{bu1}, we have
$$(x)_n \simeq (x \punt \chi)^n \simeq (x \punt \beta \punt u^{<-1>})^n \simeq [x \punt
\beta\punt  (\iota\punt \chi)_{\scriptscriptstyle D}]^n.$$ The
result follows from equivalence (\ref{(exp23)}), being
$$(x \punt \chi)^n \simeq  \sum_{k \geq 0} s(n,k) \, x^k.$$
\end{proof}

\subsection{Solving recursions}

In many special combinatorial problems, the hardest part of the solution may
be the discovery of an effective recursion. Once a recursion has been
established, Sheffer polynomials are often a simple and general tool for
finding answers in closed form. Main contributions in this respect are due to
Niederhausen \cite{Nie1,Nie2,Nie3}. Further contributions are given by Razpet
\cite{Razpet} and Di Bucchianico and Soto y Koelemeijer \cite{dibucchianico1}.

\begin{example}
\label{1} \textrm{\textrm{Suppose we are asked to solve the difference
equation
\begin{equation}
s_{n}(x+1)=s_{n}(x)+s_{n-1}(x)\label{recurr1}%
\end{equation}
under the condition $\int_{0}^{1}s_{n}(x)dx=1$ for all nonnegative integers
$n$. Equation (\ref{recurr1}) fits the Sheffer identity (\ref{(shefid1)}) if
we set $y=1,$ choose the sequence $\{p_{n}(x)\}$ such that $p_{0}%
(x)=1,p_{1}(1)=1$ and $p_{n}(1)=0$ for all $n\geq2$ and consider
the Sheffer sequence $\{n!s_{n}(x)\}$. The sequence $\{p_{n}(x)\}$
is associated to the umbra $\chi,$ so we are looking for solutions
of (\ref{recurr1}) such that
$n!s_{n}(x)\simeq(\sigma_{x}^{(\alpha,\gamma)})^{n}$ with $\gamma^{\ast}%
\equiv\chi,$ i.e. $\gamma\equiv u$ (cfr. Example \ref{(adjchi)}).
The condition $\int_{0}^{1}s_{n}(x)dx=1$ can be translated in
umbral terms by looking for an umbra $\delta$ such that
$E[s_{n}(\delta)]=1$ for all nonnegative integers $n$. Such an
umbra has g.f.
\[
\int_{0}^{1}e^{xt}dx=\frac{e^{t}-1}{t}%
\]
so that $\delta\equiv-1\boldsymbol{.}\iota,$ with $\iota$ the Bernoulli umbra.
Therefore, the umbra $\alpha$ satisfies the following identity
\[
(-1\boldsymbol{.}\alpha+-1\boldsymbol{.}\iota)\boldsymbol{.}\chi\equiv u.
\]
Due to statement \textit{i)} of Proposition \ref{prop2} , being $\beta
\boldsymbol{.}\chi\equiv u,$ we have $-1\boldsymbol{.}\alpha+-1\boldsymbol{.}%
\iota\equiv\beta$ and so $\alpha\equiv-1\boldsymbol{.}(\iota+\beta).$
Solutions of (\ref{recurr1}) are moments of the Sheffer umbra $(\iota
+\beta+x\boldsymbol{.}u)\boldsymbol{.}\chi$ normalized by $n!$.}}
\end{example}

\begin{example}
\textrm{\textrm{Suppose we are asked to solve the difference equation
\begin{equation}
s_{n}(x)=s_{n}(x-1)+s_{n-1}(x)\label{recurr2}%
\end{equation}
that satisfies the initial condition
\begin{equation}
s_{n}(1-n)=\sum_{i=0}^{n-1}s_{i}(n-2\,i)\quad\hbox{for all}\quad
n\geq1,\,\,\,s_{0}(-1)=1.\label{icr2}%
\end{equation}
If we rewrite (\ref{recurr2}) as $s_{n}(x-1)=s_{n}(x)-s_{n-1}(x),$
we note that this equation fits the Sheffer identity
(\ref{(shefid1)}) if we set $y=-1,$ choose the sequence
$\{p_{n}(x)\}$ such that $p_{0}(x)=1,$ $p_{1}(-1)=-1$ and
$p_{n}(-1)=0$ for all $n\geq2,$ and consider the Sheffer sequence
$\{n!s_{n}(x)\}$. In particular, from (\ref{(shefid)}) we have
$-1\boldsymbol{.}\gamma^{\ast}\equiv-\chi$ so $\gamma^{\ast}\equiv
-1\boldsymbol{.}-\chi$ which has g.f. $f(-1\boldsymbol{.}%
-\chi,t)=(1-t)^{-1}.$ Suppose to set $\overline{u}=-1\boldsymbol{.}-\chi$. We
have $E[\overline{u}^{n}]=n!$ for all $n\geq1$ and $x\boldsymbol{.}%
\overline{u}\equiv-x\boldsymbol{.}-\chi$ with g.f.
\[
f(x\boldsymbol{.}\overline{u},t)=(1-t)^{-x}=\sum_{n\geq0}(x)^{n}\frac{t^{n}%
}{n!}%
\]
and $(x)^{n}=x(x+1)\cdots(x+n-1)=(x+n-1)_{n}.$ In particular we
have $\gamma\equiv(\chi\boldsymbol{.}\overline{u})^{<-1>}.$
Solutions of (\ref{recurr2}) are therefore of the type
\[
s_{n}(x)\simeq\frac{\lbrack(-1\boldsymbol{.}\alpha+x\boldsymbol{.}%
u)\boldsymbol{.}\overline{u}]^{n}}{n!}.
\]
Now we need to identify $\alpha.$ As in the previous example, the moments of
such an umbra depend on the initial condition (\ref{icr2}). Observe that, if
$s_{n}(x)$ is a Sheffer sequence with associated polynomials $(x\boldsymbol{.}%
\overline{u})^{n}/n!,$ then $s_{n}(x-n)$ is a Sheffer sequence with associated
polynomials
\[
\frac{\lbrack(x-n)\boldsymbol{.}\overline{u}]^{n}}{n!}\simeq\left(
\begin{array}
[c]{c}%
x-n+n-1\\
n
\end{array}
\right)  \simeq\frac{\lbrack(x-1)\boldsymbol{.}\chi]^{n}}{n!},
\]
since $[(x-1)\boldsymbol{.}\chi]^{n}\simeq(x-1)_{n}.$ Therefore
due to Theorem  \ref{(212)} we have
\begin{equation}
s_{n}(x-n)\simeq\frac{\lbrack(-1\boldsymbol{.}\alpha+(x-1)\boldsymbol{.}%
u)\boldsymbol{.}\chi]^{n}}{n!}.\label{fin}%
\end{equation}
So the values of $s_{n}(1-n)$ in the initial condition (\ref{icr2}) give
exactly the moments of $-1\boldsymbol{.}\alpha\boldsymbol{.}\chi$ normalized
by $n!.$ Therefore, by observing that
\[
s_{n}(x)\simeq\frac{\lbrack(-1\boldsymbol{.}\alpha+x+n-1)\boldsymbol{.}%
\chi]^{n}}{n!}\simeq\sum_{k=0}^{n}\frac{(-1\boldsymbol{.}\alpha\boldsymbol{.}%
\chi)^{k}}{k!}\left(
\begin{array}
[c]{c}%
x+n-1\\
n-k
\end{array}
\right)
\]
we have
\[
s_{n}(x)=\sum_{k=0}^{n}s_{k}(1-k)\left(
\begin{array}
[c]{c}%
x+n-1\\
n-k
\end{array}
\right)  .
\]
From a computational point of view, this formula is very easy to
implement by using the recursion of the initial condition. If one
would evaluate the moments of
$-1\boldsymbol{.}\alpha\boldsymbol{.}\chi$ not by using the
recursion of the initial condition, but with a closed form, some
different considerations must be done. By using equivalence
(\ref{fin}) we have
\[
\sum_{i=0}^{n-1}s_{i}(n-2i)=\sum_{i=0}^{n-1}s_{i}(x-i)|_{x=n-i}\simeq
\sum_{i=0}^{n-1}\frac{[-1\boldsymbol{.}\alpha\boldsymbol{.}\chi
+(n-i-1)\boldsymbol{.}\chi]^{i}}{i!}%
\]
and from the initial condition we have
\begin{align}
\frac{(-1\boldsymbol{.}\alpha\boldsymbol{.}\chi)^{n}}{n!} &  \simeq\sum
_{j=0}^{n-1}\frac{(-1\boldsymbol{.}\alpha\boldsymbol{.}\chi)^{j}}{j!}%
\sum_{i=0}^{n-j-1}\frac{[(n-i-j-1)\boldsymbol{.}\chi]^{i}}{i!}\nonumber\\
&  \simeq\frac{(-1\boldsymbol{.}\alpha\boldsymbol{.}\chi+\overline{\delta
})^{\,n-1}}{(n-1)!}\label{fin1}%
\end{align}
where $\overline{\delta}$ is an umbra such that
\begin{equation}
\overline{\delta}^{\,k}\simeq k!\sum_{i=0}^{k}\frac{[(k-i)\boldsymbol{.}%
\chi]^{i}}{i!}\simeq k!\delta^{k}.\label{Fib}%
\end{equation}
Observe that
\[
\overline{\delta}^{\,k}\simeq k!\sum_{t=0}^{k}\frac{(t\boldsymbol{.}%
\chi)^{k-t}}{(k-t)!}\simeq\sum_{t=0}^{k}\left(
\begin{array}
[c]{c}%
k\\
t
\end{array}
\right)  (t\boldsymbol{.}\chi)^{k-t}\overline{u}^{t}\simeq(\overline
{u}\boldsymbol{.}\beta\boldsymbol{.}\chi_{{\scriptsize D}})^{k}%
\]
by using Corollary \ref{ccc}. The umbra $\delta$ is said to be the
boolean cumulant umbra of $\chi_{{\scriptsize D}}$ (cfr.
\cite{free}). In particular the umbra $\delta$ has moments equal
to the Fibonacci numbers, since $\overline{\delta}$ has g.f.
\[
f(\overline{\delta},t)=\frac{1}{1-t(1+t)}=\frac{1}{1-t-t^{2}}.
\]
In terms of g.f.'s, equivalence (\ref{fin1}) gives
\[
\sum_{n\geq0}\frac{(-1\boldsymbol{.}\alpha\boldsymbol{.}\chi)^{n}}{n!}%
t^{n}\simeq1+t\sum_{n\geq1}\frac{(-1\boldsymbol{.}\alpha\boldsymbol{.}%
\chi+\overline{\delta})^{\,n-1}}{(n-1)!}t^{n-1}%
\]
so that
\[
\frac{1}{f(\alpha\boldsymbol{.}\chi,t)}=1+t\frac{f(\overline{\delta}%
,t)}{f(\alpha\boldsymbol{.}\chi,t)}\Leftrightarrow f(-1\boldsymbol{.}%
\alpha\boldsymbol{.}\chi,t)=\frac{1}{1-tf(\overline{\delta},t)}%
\]
and
\[
-1\boldsymbol{.}\alpha\boldsymbol{.}\chi\equiv\overline{u}\boldsymbol{.}%
\beta\boldsymbol{.}\overline{\delta}_{{\scriptsize D}}.
\]
Therefore the solution in closed form is}}
\[
s_{n}(x)\simeq\frac{\lbrack\overline{u}\boldsymbol{.}\beta\boldsymbol{.}%
\overline{\delta}_{{\scriptsize D}}+(x+n-1)\boldsymbol{.}\chi]^{n}}{n!}.
\]

\end{example}

\begin{example}
\textrm{\textrm{Suppose we are asked to solve the difference equation
\begin{equation}
F_{n}(m)=F_{n}(m-1)+F_{n-1}(m-2)\label{recurr3}%
\end{equation}
under the condition $F_{n}(0)=1$ for all nonnegative integers $n$. Replace $m$
with $x+n+1$. Then equation (\ref{recurr3}) can be rewritten as
\begin{equation}
F_{n}(x+n+1)=F_{n}(x+n)+F_{n-1}(x+n-1).\label{recurr3bis}%
\end{equation}
Equation (\ref{recurr3bis}) fits the Sheffer identity
(\ref{(shefid1)}) if we set $y=1,$ choose the sequence
$\{p_{n}(x)\}$ such that $p_{0}(x)=1,p_{1}(1)=1$ and $p_{n}(1)=0$
for all $n\geq2$ and consider the Sheffer sequence
$\{n!F_{n}(x)\}$. As in Example \ref{1}, we are looking for
solutions of (\ref{recurr3bis}) such that
\[
F_{n}(x+n)\simeq\frac{\lbrack(-1\boldsymbol{.}\alpha+x\boldsymbol{.}%
u)\boldsymbol{.}\chi]^{n}}{n!}.
\]
Let us observe that equation (\ref{recurr3bis}), for $x=0,$ gives
the well-known recurrence relation for Fibonacci numbers so that
\[
F_{n}(n)\simeq\frac{(-1\boldsymbol{.}\alpha\boldsymbol{.}\chi)^{n}}{n!}%
\simeq\delta^{n}.
\]
Therefore we have $-1\boldsymbol{.}\alpha\boldsymbol{.}\chi\equiv
\overline{\delta},$ with $\overline{\delta}^{\,n}\simeq
n!\delta^{n}$ as given in equivalence (\ref{Fib}), and
$\alpha\equiv-1\boldsymbol{.}\overline{\delta
}\boldsymbol{.}\beta.$ Then, solutions of (\ref{recurr3bis}) are
such that
\[
F_{n}(x+n)\simeq\frac{(\overline{\delta}+x\boldsymbol{.}\chi)^{n}}{n!}.
\]
In particular
\begin{align*}
F_{n}(x+n)  &  \simeq\sum_{k=0}^{n}\frac{(x\boldsymbol{.}\chi)^{n-k}}%
{(n-k)!}\delta^{k}\simeq\sum_{k=0}^{n}\frac{(x\boldsymbol{.}\chi)^{n-k}%
}{(n-k)!}\sum_{j=0}^{k}\frac{(j\boldsymbol{.}\chi)^{k-j}}{(k-j)!}\\
&  \simeq\sum_{k=0}^{n}\sum_{j=0}^{n}\frac{(x\boldsymbol{.}\chi)^{n-k-j}%
}{(n-k-j)!}\frac{(k\boldsymbol{.}\chi)^{j}}{j!}\simeq\sum_{k=0}^{n}%
\frac{[(x+k)\boldsymbol{.}\chi]^{n-k}}{(n-k)!}\\
&  \simeq\sum_{k=0}^{n}\left(
\begin{array}
[c]{c}%
x+k\\
n-k
\end{array}
\right)
\end{align*}
by which we can verify that the initial conditions $F_{n}(0)=F_{n}(-n+n)=1$
hold.}}
\end{example}

\end{document}